\documentclass[11pt]{amsart}%{article}%
\usepackage{amsfonts,latexsym}
\usepackage{amsmath}
\usepackage{amscd}
\usepackage{amssymb}
\usepackage[pdftex]{graphicx}
\usepackage[usenames, dvipsnames, pdftex]{color}

\newcommand{\nequation}{\setcounter{equation}{0}}

\newcommand{\R}{\mathbb{R}}

\newcommand{\llangle}{\langle\!\langle}
\newcommand{\rrangle}{\rangle\!\rangle}

\newcommand{\Vir}{\text{\upshape Vir}}

\newcommand{\Diff}{\text{\upshape Diff}}
\newcommand{\Met}{\text{\upshape Met}}
\newcommand{\Vect}{\text{\upshape Vect}}
\newcommand{\Metmu}{\Met_{\mu}}
\newcommand{\Vol}{\text{\upshape Vol}}
\newcommand{\Diffmu}{\Diff_{\mu}}

\newcommand{\DiffM}{\Diff(M)}
\newcommand{\MetM}{\Met(M)}

\newcommand{\VolM}{\Vol(M)}
\newcommand{\DiffmuM}{\Diffmu(M)}

\newcommand{\Emb}{\mathcal{E}}

\newcommand{\Lie}{\mathcal{L}}
\newcommand{\ad}{\text{ad}}
\newcommand{\Ad}{\text{Ad}}
\newcommand{\grad}{\nabla}

\newcommand{\Laplacian}{\Delta}

\DeclareMathOperator{\diver}{div}
\DeclareMathOperator{\Tr}{Tr}

\DeclareMathOperator{\sgrad}{sgrad}
\DeclareMathOperator{\curl}{curl}

\newtheorem{theorem}{Theorem}[section]
\newtheorem{proposition}[theorem]{Proposition}
\newtheorem{lemma}[theorem]{Lemma}

\newtheorem{remark}[theorem]{Remark}
\newtheorem{remark-definition}[theorem]{Remark-Definition}
\newtheorem{example}[theorem]{Example}

\newtheorem{corollary}[theorem]{Corollary}

\title[Curvatures of diffeomorphism groups]{Curvatures of Sobolev metrics on diffeomorphism groups}

\author{B. Khesin}
\address{B.K.: Department of Mathematics, University of Toronto, M5S 2E4, Canada} 
\email{khesin@math.toronto.edu}

\author{J. Lenells}
\address{J.L.: Department of Mathematics, Baylor University, One Bear Place \#97328, Waco, TX 76798, USA.}
\email{Jonatan\_Lenells@baylor.edu}

\author{G. Misio\l ek}
\address{G.M.: Institute for Advanced Study, Princeton, NJ 08540, USA 
and Department of Mathematics, University of Notre Dame, IN 46556, USA}
\email{gmisiole@nd.edu}

\author{S. C. Preston}
\address{S.C.P.: Department of Mathematics, University of Colorado, CO 80309, USA}
\email{Stephen.Preston@colorado.edu}

\begin{document}

%\today 

\bigskip 
\hfill{\it To Dennis Sullivan~} 

\hfill{\it on the occasion of his 70th birthday} 
\bigskip 
\bigskip 

\begin{abstract} 
\noindent 
Many conservative partial differential equations correspond to geodesic equations 
on groups of diffeomorphisms. 
Stability of their solutions can be studied by examining sectional curvature of these groups: 
negative curvature in all sections implies exponential growth of perturbations and 
hence suggests instability, 
while positive curvature suggests stability. 
In the first part of the paper 
we survey what we currently know about 
the curvature-stability relation in this context 
and provide detailed calculations for several equations of continuum mechanics 
associated to Sobolev $H^0$ and $H^1$ energies. 
In the second part 
we prove that in most cases (with some notable exceptions) 
the sectional curvature assumes both signs. 
\end{abstract}

\maketitle



\tableofcontents

%%%%%%%%%%%%%%%%%%%%%%%%%%%%%%%%%%%%%%%%%%
%%%%%%%%%%%%%%%%%%%%%%%%%%%%%%%%%%%%%%%%%%%%%%
\section{Introduction} 
\nequation

The idea that stability of a dynamical system can be investigated 
using tools of Riemannian geometry goes back to Hadamard~\cite{hadamard} 
who studied the free motion of a particle on a surface of constant negative curvature. 
A similar approach was considered by Synge \cite{synge}. 
Perhaps the most influential example is due to Arnold \cite{A} who showed that 
fluid motions can be viewed as geodesics in the infinite-dimensional group of 
volume-preserving diffeomorphisms. This led him to examine curvature of 
the diffeomorphism group and derive a number of results on stability of ideal fluids. 
Roughly speaking, on a finite-dimensional manifold negative sectional curvature 
is related to instability,
and positive sectional curvature is related to stability of the corresponding geodesic flow;
the Rauch comparison theorem makes this comparison rigorous. 
Since the work of Arnold other partial differential equations 
have been interpreted as geodesic equations in infinite-dimensional spaces 
and, as with ideal fluids, calculating sectional curvatures in these cases 
has become a matter of broader interest. 

In this paper we examine the sign of the sectional curvature of certain metrics 
on infinite dimensional manifolds
(which are associated with several well-known equations of mathematical physics)
and its relevance in the stability analysis of the associated initial value problems. 

Our main interest is in those equations that arise from right-invariant metrics 
on the group of (smooth) diffeomorphisms $\Diff(M)$,
or its subgroup of volume-preserving diffeomorphisms (volumorphisms) $\Diffmu(M)$,
of a compact $n$-dimensional manifold $M$ without boundary. 
Both spaces can be completed to Hilbert manifolds $\Diff^s(M)$ and $\Diff^s_{\mu}(M)$ 
modelled on Sobolev spaces of $H^s$ vector fields and divergence-free vector fields, respectively, 
with $s>n/2 + 1$. 
We assume that  $M$ has a Riemannian metric $\langle \cdot, \cdot \rangle$ 
with volume form $\mu$. 
We equip the groups with right-invariant Sobolev metrics 
such that, on the tangent space at a diffeomorphism $\eta$, we have 
\begin{equation} \label{rightinvmetric}
\llangle u\circ\eta, v\circ\eta\rrangle_{\eta} 
= 
\int_M \Big( a\langle u, v\rangle + b\langle \delta u^{\flat}, \delta u^{\flat}\rangle 
+ c\langle du^{\flat}, du^{\flat}\rangle\Big) \, d\mu,
%\int_M \langle u, \Lambda v\rangle \, d\mu 
\end{equation} 
for any vector fields $u, v$ on $M$. 
%and $\Lambda$ is a positive-definite symmetric operator on 
%We will focus on the case of Sobolev $H^1$ products with 
%
%\begin{equation}\label{abclambda}
%\Lambda u = a u + ( b d \delta u^{\flat}+ c \delta du^{\flat} )^{\sharp},
%\end{equation}
%
Here $d$ is the exterior derivative, $\delta = \pm {*} d*$ 
%\footnote{\color{Green} What happened to the factor $(-1)^{nk +1}$ in the definition
%$$\delta \beta = (-1)^{nk +1}{*}d{*}\beta, 
%\qquad 
%\beta \in \Omega^{k+1}(M)?$$
%Why do we have just $-1$?
%It should be there. I think it might have been originally because we just did it on 1-forms or something. I think the sign is unimportant so I changed it to just +/-.}
is its (formal) adjoint, $a,b,c$ are non-negative constants and 
$\flat$ and $\sharp$ denote the standard ``musical isomorphisms'' of the metric 
corresponding to lowering and raising of indices. 
Formula \eqref{rightinvmetric} simplifies in dimensions $n \leq 3$ 
where we have\footnote{Recall that $\curl{u} = *du^\flat$ and $\sgrad{f}= -(\delta{*}f)^\sharp$ in two dimensions, while $\curl{u} = (*du^\flat)^\sharp$ in three dimensions.}
\begin{equation}\label{rightinvmetricspecialcases} 
\llangle u\circ\eta, v\circ\eta\rrangle_\eta 
= 
\begin{cases} \displaystyle
\int_M ( auv + bu_x v_x ) \, dx & n=1,  \\  \displaystyle 
\int_M \big( a\langle u, v\rangle + b\diver{u}\cdot \diver{v} + c \curl{u} \cdot \curl{v} \big) \, dA & n=2, 
\\ \displaystyle 
\int_M \big(a\langle u, v\rangle + b\diver{u}\cdot \diver{v} + c\langle \curl{u}, \curl{v}\rangle\big) \, dV & n=3.
\end{cases}
\end{equation}
%
%\begin{alignat}{2}  \nonumber 
%n&=1 \quad & \quad \llangle u\circ\eta, 
%\Lambda u &= a u - b u''   %\label{1dlambda} 
%\\  \nonumber
%n&=2 \quad & \Lambda u &= a u - b \grad \diver{u} + c \sgrad \curl{u}   %\label{2dlambda}  
%\\ \nonumber
%n&=3 \quad & \Lambda u &= a u - b\grad \diver{u} + c \curl(\curl{u}).   %\label{3dlambda}
%\end{alignat}
%

Using the metric \eqref{rightinvmetric} one derives a number of PDE 
that are of interest in continuum mechanics and geometry. 
For example, in the one-dimensional case one obtains 
Burgers' equation, the Camassa-Holm equation, and the Hunter-Saxton equation 
as geodesic equations for appropriate choices of $a$, $b$, $c$. 
In higher dimensions one gets the EPDiff equation and the so-called template-matching equation. 
Furthermore, projection onto $\Diffmu(M)$ yields the usual Euler equations of hydrodynamics 
and the Lagrangian-averaged Euler-$\alpha$ equation. 
We refer to the paper \cite{klmp} and the references therein for more details. 

The goal of this paper is two-fold. In the first part we describe some aspects of 
Riemannian geometry of infinite-dimensional manifolds which are relevant to 
the analysis of partial differential equations of mathematical physics. 
In particular, we review the framework of the Euler-Arnold equations 
on Lie groups equipped with right-invariant metrics and explain the role 
played by sectional curvature in the study of (Lagrangian) stability. 

In the second part we present new results on the sign of the sectional curvature 
for the metric \eqref{rightinvmetric} on $\Diff(M)$ and $\Diffmu(M)$
for different choices of the parameters $a$, $b$ and $c$. 
In order to simplify calculations, we will take $M$ to be either the circle $S^1$ 
or the flat torus $\mathbb{T}^n$. 
We show that in most cases the sectional curvature assumes both signs. 
Two notable exceptions are: 
$b=c=0$ and $n=1$, in which case the sectional curvature of 
$\Diff(S^1)$ turns out to be non-negative; 
and 
$a=c=0$ and $n \geq 1$, in which case the sectional curvature of 
the quotient space $\Diff(M)/\Diffmu(M)$ is strictly positive for any compact $M$. 
The latter case is studied in detail in the  paper \cite{klmp}. 
For the $H^1$ metric on $\Diff(S^1)$ we obtain 
a simple curvature expression and explain how it can be viewed as 
the Gauss-Codazzi formula for an isometric embedding of the group of circle diffeomorphisms 
in a larger space.

Besides the fact that right-invariant $H^1$ metrics \eqref{rightinvmetric} arise as 
Lagrangians of many PDE of continuum mechanics, 
our motivation to study them is also purely geometric: they are the natural metrics 
induced on orbits of pullback actions on spaces of tensor fields. 
For example, the canonical $L^2$ metric on the space of all Riemannian metrics 
on a compact manifold $M$ 
induces a metric of the type \eqref{rightinvmetric} on $\Diff(M)$ (viewed as an orbit of 
any particular metric under the pullback action)
whose geometry was studied in \cite{klmp}. 
%with operator \eqref{abclambda} 
%for which $a=0$, $b=4$, and $c=2$, see Remark 8.3 of \cite{klmp}.
%with a particular choice of constants $a,b$ and $c$, see \cite{klmp}.

%%%%%%%%%%%%%%%%%%%%%%%%%%%%
%%%%%%%%%%%%%%%%%%%%%%%%%%%%%%%%%%%%%%%%%%%

\section{Preliminaries: metrics and geodesic equations} 
\label{sec:Back} 
\nequation 

In this section we review some well-known examples of infinite-dimensional Riemannian manifolds 
and their geodesic equations. 
Most interesting from our point of view are infinite-dimensional Lie groups, 
especially diffeomorphism groups equipped with right-invariant metrics. 
We describe examples of such groups in Section \ref{groupsubsection}. 
In Section \ref{othermfdsubsection} we discuss other situations of interest 
when the manifolds are not Lie groups or the metrics are not right-invariant. 

%%%%%%%
\subsection{The Euler-Arnold equations on Lie groups: examples} 
\label{groupsubsection}

We begin by describing the general Lie-theoretic setup of Arnold \cite{A}. 

Consider an infinite-dimensional Lie group $G$ equipped with 
a smooth right-invariant (weak) Riemannian metric determined by 
an inner product $\llangle \cdot, \cdot \rrangle$ defined on the tangent space $T_eG$ 
at the identity element. 
A geodesic in the group starting from $e$ in the direction $u_0$ can be obtained from 
the solution $u(t)$ of the Cauchy problem for the associated Euler-Arnold equation 
on $T_e G$, namely 
\begin{equation} \label{utBuu} 
\frac{du}{dt} = - \ad^\ast_uu, 
\qquad 
u(0)=u_0 
\end{equation} 
where $u \to \ad^\ast_v u$ is 
the adjoint of the linear operator $u \to \ad_vu = -[v,u]$ 
with respect to the inner product on $T_eG$, that is\footnote{Here $[\cdot, \cdot]$ denotes 
the commutator on $T_eG$ induced by the Lie bracket of right-invariant vector fields on $G$, i.e. 
$[u, v] = [X,Y]_e$ where $X,Y$ are the right-invariant vector fields determined by $X_e = u$ and $Y_e = v$.}
\begin{equation} \label{opB} 
\llangle \ad^\ast_v u, w\rrangle 
= 
-\llangle u, [v,w]\rrangle, 
\qquad 
u, v, w \in T_e G. 
\end{equation} 
The geodesic is now obtained by solving the flow equation 
\begin{equation} \label{floweq} 
\frac{d\eta}{dt} = DR_{\eta(t)} u(t), 
\qquad 
\eta(0)= e 
\end{equation} 
where $\xi \mapsto R_\eta(\xi)$ denotes the right-translation in $G$ by $\eta$. 

In the case of ideal hydrodynamics $G$ is the group of volumorphisms 
(volume-preserving diffeomorphisms) of a manifold $M$ 
$$
\DiffmuM = \{ \eta\in \Diff(M) \, \vert \, \eta^*\mu = \mu\},
$$
with the right-invariant metric given at the identity by the $L^2$ inner product, 
i.e. by setting $b=c=0$ in \eqref{rightinvmetric}. 
The resulting Euler-Arnold equations \eqref{utBuu} are the familiar Euler equations of 
incompressible fluids in $M$ 
\begin{equation}\label{E}  
u_t + \nabla_u u = - \nabla p, 
\quad 
\diver{u} = 0,
\end{equation} 
frequently written in the form 
\begin{equation}\label{vorticityeuler}
\omega_t + \Lie_u\omega = 0, 
\end{equation} 
where $\omega = \curl{u}$ is the vorticity.\footnote{In two dimensions 
$\omega=\curl{u}$ is a function and $\Lie_u\omega = \langle u, \grad \omega\rangle$; 
in three dimensions $\omega=\curl{u}$ is a vector field and $\Lie_u\omega = [u,\omega]$.}

If we set $a=1$ and $c=\alpha^2$ with $b=0$ in \eqref{rightinvmetric} 
on the volumorphism group 
then the corresponding Euler-Arnold equation \eqref{utBuu} 
is called the Lagrangian-averaged Euler equation \cite{HMR, shkoller}. 
It is more complicated than \eqref{E}, but its ``vorticity'' given by 
$\omega = \curl{u} - \alpha^2 \Laplacian \curl{u}$ 
satisfies the same equation as \eqref{vorticityeuler}. 
The analysis of this equation presents similar difficulties as \eqref{E} or 
\eqref{vorticityeuler}; see e.g., \cite{houli}.

\begin{remark}\upshape
The volumorphism group is one of the three ``classical'' diffeomorphism groups. Another is the symplectomorphism group $\Diff_{\omega}(M)$, 
consisting of diffeomorphisms preserving a symplectic form $\omega$ of an even-dimensional manifold $M$,
and the third is the contactomorphism group $\Diff_{\alpha}(M)$ 
consisting of diffeomorphisms $\eta$ such that $\eta^*\alpha = F\alpha$, 
where $\alpha$ is a contact form and $F$ is a nowhere-zero function 
on an odd-dimensional manifold $M$. 
Geodesic equations of the right-invariant $L^2$ metric on these groups have been studied 
in \cite{ebinsymplectic, khesinsymplectic, ebinpreston}. 
In two dimensions the geodesic equation on the symplectomorphism group 
reduces to the Euler equation \eqref{E}, and in one dimension 
the geodesic equation for contactomorphisms reduces to the Camassa-Holm equation. 
Simpler equations arise when the $L^2$ metric is restricted to 
the subgroup of Hamiltonian diffeomorphisms or the subgroup of strict contactomorphisms; 
see Smolentsev~\cite{sm} for a review of their properties. 

It is worth pointing out that the subgroup of Hamiltonian diffeomorphisms 
carries a bi-invariant metric given at the identity by 
\begin{equation}\label{biinvariant}
\llangle \sgrad f, \sgrad g\rrangle = \int_M fg \, d\mu,
\end{equation}
where $f$ and $g$ are assumed to have mean zero.
For such metrics we have 
$\llangle u, \ad_vw\rrangle + \llangle \ad_vu, w\rrangle = 0$ 
whenever $u, v$ and $w$ are in the Lie algebra. 
It follows that $\ad_u^\ast u=0$, and hence the Euler equation \eqref{utBuu} 
reduces to $du/dt=0$. 
Geodesics are easy to find as they are one-parameter subgroups: 
simply fix a velocity field and compute the flow. 
\end{remark}

%If $M$ is an even-dimensional manifold with a symplectic form $\omega$, then we can consider a right-invariant $L^2$ metric on the group $\Diff_{\omega}(M) = \{\eta \in \Diff(M) \, \vert \, \eta^*\omega = \omega\}$ of symplectomorphisms, or on the subgroup of Hamiltonian diffeomorphisms (these coincide if the first cohomology $H^1(M, \mathbb{R})$ is trivial). 
%The geodesic equation on these groups has been studied by Ebin~\cite{ebinsymplectic} in the former case and by Smolentsev~\cite{smolentsev_symplecto} in the latter case; in two dimensions it reduces to the Euler equation \eqref{E} and has similar analytic properties (such as global existence) in any dimension.

%If $M$ is an odd-dimensional manifold with a contact form $\alpha$, we can similarly consider a right-invariant $L^2$ metric on the contactomorphism group $$
%\Diff_{\alpha}(M) = \{ \eta\in \Diff(M) \, \vert \, \eta^*\alpha = F\alpha \text{ for some $F\in C^{\infty}(M,\mathbb{R}^+)$}\}.$$
%In one dimension this equation reduces to the Camassa-Holm equation; it is studied in \cite{ebinpreston}. The analogous equation for the strict contactomorphism group (or quantomorphism group) consisting of those diffeomorphisms strictly preserving $\alpha$ was analyzed by Smolentsev~\cite{smolentsevclassical}.

The group of circle diffeomorphisms $\Diff(S^1)$ has been a rich source of examples. 
In this case we can set $c=0$ in the formula \eqref{rightinvmetric} since in one dimension 
$du^{\flat}=0$ for any vector field $u$. 
Two much-studied Euler-Arnold equations that arise here are 
the periodic (inviscid) Burgers equation 
\begin{equation}\label{B} 
u_t + 3uu_x = 0
\end{equation} 
associated with the $L^2$ inner product (with $b=0$) 
and the periodic Camassa-Holm equation 
\begin{equation}\label{CH} 
u_t - u_{txx} + 3uu_x - 2u_x u_{xx} - uu_{xxx} = 0 
\end{equation} 
obtained from the $H^1$ product (with $a=b=1$) \cite{ak, misiolekcamassa}. 

Interesting examples also arise on the Bott-Virasoro group $\Vir(S^1)$, 
the universal  central extension of $\Diff(S^1)$,
with a group law defined by 
\begin{equation} \label{centralextensiongrouplaw}
(\eta, \alpha)\circ (\xi, \beta) 
= 
\left(\eta\circ\xi, \alpha + \beta +
\frac{1}{2}\int_{S^1}\log{\partial_x(\eta\circ\xi)}\,d\log{\partial_x\xi}\right) 
\end{equation}
where $\eta,\xi\in \Diff(S^1)$ and $\alpha, \beta \in \R$.
The right-invariant metric on the Bott-Virasoro group given at the identity by 
the $L^2$ inner product 
$$
\llangle (u, a),(v, b) \rrangle_{L^2} = \int_{S^1} u v \, dx + ab,
$$
yields as its Euler-Arnold equation \eqref{utBuu} the Korteweg-de Vries equation 
\cite{ak, OK}
\begin{equation}\label{kdv} 
u_t + 3 u u_{x} + a u_{xxx} = 0, \quad a = \mathrm{const}.
\end{equation} 
If we use the $H^1$ inner product instead 
$$ 
\llangle (u, a), (v, b)\rrangle_{H^1} = \int_{S^1} (uv + u_x v_x) \, dx + ab,
$$
then the corresponding Euler-Arnold equation is the Camassa-Holm equation 
with drift~\cite{misiolekcamassa} 
$$ 
u_t - u_{txx} + 3uu_x - 2u_x u_{xx} + (\kappa - u) u_{xxx} = 0, \quad \kappa = \text{const}.
$$

%Similarly, the quotient space $\Diff(S^1)/\Rot{S^1}$ is a configuration space of 
%the periodic Hunter-Saxton equation 
%
%\begin{equation} \tag{HS} 
%u_{txx} + 2u_x u_{xx} + uu_{xxx} =0
%\end{equation} 
% 
%which arises from the (homogeneous) $\dot{H}^1$ product with $a=0$
%on the space of mean-zero periodic functions. 

Equations \eqref{B} and \eqref{CH} have higher-dimensional analogues. 
The $n$-dimensional version of the Burgers equation \eqref{B} arising from a right-invariant $L^2$ metric on $\Diff(M)$ 
is the so-called template-matching equation~\cite{mumford} 
$$ 
u_t + \nabla_uu + (\nabla u)^{\dagger}(u) + (\diver{u})u = 0 
$$
while the $n$-dimensional version of the Camassa-Holm equation \eqref{CH} on $\Diff(M)$ is the EPDiff equation~\cite{HMR}.

%The Hunter-Saxton equation also has a remarkable higher dimensional generalization, discussed in \cite{klmp}. 
%The quotient $G=\DiffM/\DiffmuM$ can be viewed as the space of 
%(smooth) positive densities on $M$ of total mass $\mu(M)$ and 
%the right-invariant metric obtained from the Sobolev $H^1$ inner product 
%by setting $a=c=0$ in \eqref{abclambda} 
%coincides with the so-called Fisher-Rao (information) metric in geometric statistics.  
%The corresponding Euler-Arnold equation has the form 
%
%\begin{equation*}
%\rho_t + \nabla_u \rho + \frac{1}{2} \rho^2 + \frac{1}{\mu(M)}\int_M \rho^2 d\mu = 0, 
%\quad 
%\rho = \diver{u}. 
%\end{equation*} 
% 
%We refer to \cite{klmp} for the derivation and detailed study of 
%this equation as well as the attendant geometry. 

%Finally, if $G$ is the Virasoro group $\Vir(S^1)=\Diff(S^1)\times\mathbb{R}$ (i.e. the universal central extension of $\Diff(S^1)$) 
%with the right-invariant $L^2$ metric the Euler-Arnold equation \eqref{utBuu} 
%becomes the Korteweg-de Vries equation 
%
%\begin{equation} \tag{KdV} 
%u_t + 3uu_x + u_{xxx} = 0.
%\end{equation} 
%

Other examples of geodesic equations for right-invariant metrics on infinite-dimensional groups 
include 
\begin{itemize}
\item 
the equation of passive scalar motion 
(on the semidirect product of $\Diffmu(M)$ with $C^{\infty}(M)$); 
\item 
the equation of $3D$ magnetohydrodynamics 
(on the semi-direct product of $\Diffmu(M)$ with divergence-free vector fields) 
$$ 
u_t + \nabla_uu = -\grad p + \nabla_BB, \qquad B_t + [u,B] = 0, \qquad \diver{u}=\diver{B}=0; 
$$
\item 
the $\mu$-CH equation, corresponding to an $H^1$-type metric on $\Diff(S^1)$, 
see Remark \ref{mu-CH}; 
\item 
the quasigeostrophic equation in $\beta$-plane approximation on $\mathbb{T}^2$ 
on the central extension of the group of Hamiltonian diffeomorphisms 
$$ 
\omega_t + \{\psi, \omega\} = -\beta \partial_x\psi, \qquad \omega = \Laplacian \psi; 
$$
\end{itemize}
as well as the Boussinesq approximation to stratified fluids, equations for charged fluids 
and fluids in Yang-Mills fields, the Landau-Lifschitz equations 
as well as 
various $2$-component generalizations of the one-dimensional equations mentioned above;
see e.g., \cite{KW, vizman}.

%%%%%%%%%%%%%%%%%%%%%%%%%%%%%%%%%%%

\subsection{Further examples of infinite-dimensional geodesic equations} 
\label{othermfdsubsection}

%Conservative second-order evolution PDE are often be derived by minimizing a Lagrangian, and if the Lagrangian happens to be quadratic in the velocities, we formally have a geodesic equation on an infinite-dimensional manifold. If this manifold happens to have a group structure under which the Riemannian metric is one-sided invariant, we are in the situation of the previous subsection, and the geodesic equations decouples into two first-order partial differential equations. However there are some infinite-dimensional manifolds in which there is no group structure (or where there is a group structure not respected by the metric) where we still care about the geodesics. In this subsection we will survey the formal properties of such geodesics, reserving the more rigorous considerations for the next subsection.

%%%%%
\subsubsection{Spaces of curves} 
Let $\Omega M$ be the (free) loop space over a compact Riemannian manifold $M$ 
whose points are smooth maps from $S^1$ to $M$.\footnote{If $M$ is a Lie group then 
$\Omega M$ becomes a loop group under pointwise multiplication.} 
The tangent space to $\Omega M$ at a point $\gamma$ consists of vector fields in $M$ 
along $\gamma$, i.e., maps $s \to V(s)\in T_{\gamma(s)}M$. 

Two metrics on $\Omega M$ have been of particular interest. 
The first is the weak Riemannian $L^2$ metric given at $\gamma \in \Omega M$ by 
\begin{equation} \label{L2loopspace}
\llangle U,V\rrangle_{L^2, \gamma} 
= 
\int_{S^1} \langle U(s), V(s)\rangle_{\gamma(s)} \, ds. 
\end{equation}
Its geodesics correspond to geodesics on the underlying manifold: if $\eta(0)=\gamma$ 
and $\dot{\eta}(0)=V$, then $\eta(t)(s) = \exp_{\gamma(s)}(tV(s))$. 

The other is the Sobolev $H^1$ metric 
\begin{equation}\label{H1loopspace}
\llangle U, V\rrangle_{H^1 , \gamma} 
= 
\int_{S^1} \left(
\langle U(s), V(s)\rangle_{\gamma(s)} 
+ 
\Big\langle \frac{DU}{ds}, \frac{DV}{ds}\Big\rangle_{\gamma(s)} 
\right) ds,
\end{equation}
where $D/ds$ denotes the covariant derivative along $\gamma$ in $M$. 
The metric \eqref{H1loopspace} is in fact the more natural of the two 
and, in particular, turns the set $\overline{\Omega M}^{_{H^1}}$ 
consisting of all $H^1$ loops in $M$ into a complete Hilbert Riemannian manifold. 

The set of simple closed curves in $\mathbb{R}^2$ can also be regarded 
as a version of the loop space. 
Such curves may be viewed as boundaries of planar ``shapes'' and 
finding a suitable notion of distance between ``shapes'' 
has been of interest in applications to pattern theory. 
For this purpose, however, parameterizations are irrelevant; thus
 it is useful to pass to the quotient by the diffeomorphism group of $S^1$. 
Geodesics on the quotient of the metrics induced by 
\eqref{L2loopspace} and \eqref{H1loopspace} 
were studied by Michor and Mumford~\cite{MM}. 
Another approach is to consider 
the subspace consisting of those curves parameterized by arc length (or its multiple) 
with the induced $L^2$ or $H^1$ metric. 
The $L^2$ geodesics on this subspace are solutions of a wave-like equation 
\begin{equation}\label{whip}
 \eta_{tt} = \partial_s(\sigma \eta_s), 
 \qquad 
 \sigma_{ss}-\lvert \eta_{ss}\rvert^2 \sigma 
 = 
 -\lvert \eta_{st}\rvert^2, \qquad \lvert \eta_s\rvert\equiv 1,
\end{equation}
which describes an inextensible string (the asymptotic limit of a string with a very strong tension);
see \cite{prestonwhip} for a geometric discussion of this equation.

%%%%%%
\subsubsection{Homogeneous spaces} 
Degenerate metrics on diffeomorphism groups also lead to natural geometries 
on their quotient spaces. 
For example, 
if $u$ is a vector field on $S^1$ and $\eta \in \Diff(S^1)$, then setting 
$$ 
\llangle u\circ\eta, v\circ\eta \rrangle_{\dot{H}^1} 
= 
\int_{S^1} u_x v_x \, dx 
$$
one obtains an invariant degenerate Sobolev $\dot{H}^1$ metric which is 
a limiting case of the $H^1$ metric \eqref{rightinvmetricspecialcases} 
when $a\to 0$ or $b\to\infty$. 
It becomes a weak Riemannian metric on the homogeneous space $\Diff(S^1)/S^1$. 
The corresponding Euler-Arnold equation is the Hunter-Saxton equation~\cite{K-M} 
\begin{equation} \label{HS1d}
u_{txx} + 2u_xu_{xx} + u u_{xxx} = 0.
\end{equation}
Passing to the quotient space is geometrically appealing since the manifold turns out 
to be isometric to a subset of the round Hilbert sphere~\cite{L}. 

More generally, the same construction applies on any Riemannian manifold $M$ 
using the right-invariant degenerate metric 
$$ 
\llangle u\circ\eta, v\circ\eta\rrangle_{\dot{H}^1} 
= 
\int_M \diver{u}\cdot \diver{v} \, d\mu,
$$
on the quotient space $\Diff(M)/\Diffmu(M)$. Its Euler-Arnold equation is 
a higher-dimensional analogue of \eqref{HS1d} given by
\begin{equation}\label{HSmultid}
\grad \diver{u}_t + \diver{u} \grad \diver{u} + \grad \langle u, \grad \diver{u}\rangle = 0,
\end{equation}
and one can establish a similar isometry with the round sphere in a Hilbert space. 
The induced Riemannian distance turns out to be a spherical analogue of 
the Hellinger metric in probability theory; see \cite{klmp} for details.

%%%%%
\subsubsection{Spaces of maps and non-invariant metrics} 
More generally, 
given a Riemannian manifold $M$ and a compact manifold $N$ with a volume form $\nu$ 
(and possibly with boundary) 
consider the space $C^{\infty}(N,M)$ of smooth maps from $N$ into $M$. 
%(This can be viewed as a generalization of the loop space, where $N=S^1$, as well as 
%the diffeomorphism group, where $N=M$.) 
%For example, when $\Omega$ is the closure of an open subset of $\mathbb{R}^n$ and $M=\mathbb{R}^n$, this is the natural configuration space for the free boundary problem of a compressible fluid. 
On each tangent space at $f \in C^\infty(N,M)$ we can define an $L^2$ metric by 
\begin{equation} \label{noninvariantL2}
\llangle U, V \rrangle_f 
= 
\int_M \langle U(x), V(x) \rangle_{f(x)} \, d\nu(x),
\end{equation}
and, as in the case of the loop space and \eqref{L2loopspace}, show that 
its geodesics come directly from geodesics on $M$. 
%We could also define an $H^1$ metric as in \eqref{H1loopspace}, although it is somewhat more complicated and less natural. 

The group of smooth diffeomorphisms $\Diff(M)$ is an open subset of 
the Frechet manifold $C^{\infty}(M, M)$ 
so that we can likewise put the metric \eqref{noninvariantL2} on it. 
Note that this metric is \emph{not} right-invariant;
nevertheless the corresponding geodesic equation can be 
rewritten on the tangent space to the identity where it becomes 
the multidimensional inviscid Burgers (or pressureless compressible Euler) equation 
\begin{equation}\label{pressurelesseuler}
\frac{\partial u}{\partial t} + \nabla_uu = 0.
\end{equation}
Since geodesics in $M$ starting from two nearby points will inevitably cross 
(at which time the geodesic in $\Diff(M)$ must exit the diffeomorphism group) 
solutions of the pressureless Euler equation solutions in general will blow up in finite time. 
Physically, this corresponds to the emergence of a shock wave 
leading to collisions of the fluid (or gas) particles. 
Nonetheless, the geodesic remains in $C^{\infty}(M, M)$ for all time.

\begin{remark} \upshape 
The equations of incompressible fluids with boundary can be viewed formally 
as geodesic equations on the space
$\Emb_{\mu}(\Omega, \mathbb{R}^n)$
 of volume-preserving embeddings  
of the closure $\bar{\Omega}$ of an open subset of $\mathbb{R}^n$ into $\mathbb{R}^n$. 
These equations were studied geometrically by Ebin~\cite{ebinfreeboundary} 
and shown to be 
identical to the standard equations of incompressible fluid mechanics 
except for the fact that the boundary condition for the pressure is 
a Dirichlet rather than a Neumann condition.
\end{remark} 

\subsubsection{Spaces of metrics}
Another geometrically interesting space is 
the space $\Met(M)$ of all Riemannian metrics on a compact manifold $M$. 
For any metric $g$ on $M$ and any $\eta \in \DiffM$ we define the pullback metric $\eta^*g$. 
In \cite{ebinmetrics} Ebin studied the Riemannian metric on $\Met(M)$ which is 
right-invariant under the pullback action. 
Given $g\in \Met(M)$ and tangent vectors $A$ and $B$ 
(smooth tensor fields of symmetric bilinear forms) the metric is defined by 
$$
\llangle A, B\rrangle_g = \int_M \Tr_g(AB) \, d\mu_g,
$$
where $\Tr_g$ is the trace with respect to $g$ and $\mu_g$ is the Riemannian volume form. 
The curvature and geodesics of this metric were computed explicitly by 
Freed and Groisser~\cite{fg}: 
sectional curvature is non-positive and geodesics generally exist only for finite time 
(until the metric becomes degenerate). 
The diffeomorphism group embeds in $\MetM$ as an orbit of a generic $g$ 
(i.e. with no non-trivial isometries) and if $g$ is Einstein then the induced metric on $\Diff(M)$ 
is a special case of \eqref{rightinvmetric}. 
We refer to Clarke~\cite{clarke} for recent results on the distance and diameter of this space.

Similarly, one can endow the space of all volume forms $\Vol(M)$ on $M$ 
with a natural right-invariant metric given for $n$-form fields $\alpha$ and $\beta$ 
tangent to $\mu$ by 
$$ 
\llangle \alpha, \beta\rrangle_{\nu} 
= 
\int_M \frac{\alpha}{\nu} \, \frac{\beta}{\nu} \, d\nu.
$$ 
Although this metric is flat, it too is not geodesically complete in general. 
Orbits of the diffeomorphism group in $\Vol(M)$ are the homogeneous spaces of 
densities $\Diff(M)/\Diffmu(M)$ of constant positive curvature; 
we refer to \cite{klmp}. 

The map $g\mapsto \mu_g$ where $\mu_g$ is the Riemannian volume of $g$ 
is a submersion which becomes a Riemannian submersion after suitable rescaling 
of the metric on $\VolM$. 
Its fibers are the spaces of metrics $\Metmu(M)$ with the same volume form $\mu$. 
These fibers are globally symmetric (with negative curvature and indefinitely extendable geodesics) 
in the induced metric from $\MetM$. 
The natural action on $\Metmu(M)$ is pullback by volumorphisms,
and if $g$ has no nontrivial isometries, then 
the orbits of $\Diffmu(M)$ are embedded submanifolds with right-invariant metrics.

%%%%%%%%%%%%%%%%%%%%%%%%%%%%%%%%
%%%%%%%%%%%%%%%%%%%%%%%%%%%%%%%%%%%%%

\section{Global aspects of infinite-dimensional Riemannian geometry} 
\label{rigoroussubsection} 
\nequation 

The obstacles that arise in the study of global Riemannian geometry of 
infinite dimensions manifolds are well known. They are mostly caused by 
the lack of local compactness or the fact that the topology generated by the metric 
may be weaker than the manifold topology. 
As a result some of the finite-dimensional techniques are not available 
or are of limited use. For example,
the Riemannian exponential map may not be defined on the whole tangent bundle 
or even be smooth, 
conjugate points may cluster along finite geodesic segments or have infinite multiplicity, etc. 
In this section we illustrate some of these situations with a few familiar examples. 

%%%%%%%%%%%
\subsection{Degenerate distance functions}\label{degeneratedistance}
The distance between two points in a weak Riemannian Hilbert manifold can be defined 
as in finite dimensions, i.e., as the infimum of lengths of piecewise smooth curves joining them. 
It is easy to prove that it satisfies all the axioms of a metric space 
except for nondegeneracy which typically requires some additional assumptions. 

As an example, consider the right-invariant $L^2$ metric on $\Diff(M)$ defined by 
$$ 
\llangle u\circ\eta, v\circ\eta\rrangle_{\eta} = \int_M \langle u,v\rangle \, d\mu, 
$$ 
which corresponds to the case $b=c=0$ in \eqref{rightinvmetric}. 
In \cite{MMvanishingdiffeo} it is shown that the geodesic distance is identically zero,  
(i.e., between any two diffeomorphisms there are curves of arbitrarily short length). 
This is essentially related to the lack of control over the Jacobian. 
The same phenomenon also occurs
for the right-invariant $L^2$ metric on the Bott-Virasoro group \cite{BBHM}, 
for the $L^2$ metric on the ``shape space'' of curves modulo reparameterizations \cite{MM},
and 
for the bi-invariant Hofer-type $L^2$ metric \eqref{biinvariant} on the Hamiltonian diffeomorphisms \cite{EP}. 

On the other hand, we obtain nondegenerate Riemannian distances for  
the $L^2$ metric on $\DiffmuM$ \cite{Ebin-Marsden}, 
the right-invariant metric on $\DiffM$ corresponding to $c=0$ in \eqref{rightinvmetric} \cite{MP}, the 
space of maps $C^{\infty}(N, M)$ in the $L^2$ metric \eqref{noninvariantL2}, 
the space of arc-length parameterized curves \cite{prestonwhip}, and the space
$\Diff(M)/\Diffmu(M)$ of densities \cite{klmp}.

%This phenomenon is essentially a consequence of the metric being not quite strong enough to handle the group structure. For the non-invariant $L^2$ metric \eqref{noninvariantL2}, the geodesics on $C^{\infty}(N, M)$ are directly related to the geodesics on $M$, and thus so is the geodesic distance on $C^{\infty}(N, M)$: we have 
%$$ D(\eta, \xi) = \sqrt{\int_N d(\eta(p), \xi(p))^2 \, d\nu(x)},$$
%where $d$ is the geodesic distance in $M$. Hence the only way to have $D(\eta, \xi)=0$ is if $\eta=\xi$. Similarly on the diffeomorphism group $\Diff(M)$ with the non-invariant metric, the geodesic distance is nondegenerate.

\subsection{Completeness and minimizing geodesics} 
Even if the geodesic distance is nondegenerate, thus providing a genuine metric space structure 
on a space of maps, this metric space may not be complete. 
For example, the completion of $C^\infty(N,M)$ in the non-invariant $L^2$ metric \eqref{noninvariantL2} 
consists of measurable maps from $N$ to $M$ which may not even be continuous. 
The same phenomenon occurs for the group of volumorphisms $\Diffmu(M)$ with the $L^2$ metric 
if $M$ is a three-dimensional manifold---in this case  the completion in the Riemannian distance 
is the space of all measure-preserving maps; see \cite{shn}.\footnote{If $M$ is two-dimensional, the completion of $\Diffmu(M)$ in the $L^2$ metric is unknown.}

In finite dimensions, completeness of a Riemannian manifold $M$ as a metric space 
is equivalent to geodesic completeness, i.e., extendability of geodesics for all time, 
which in turn implies that any two points in $M$ can be joined by a minimal geodesic. 
The proof of this result, the Hopf-Rinow theorem, 
relies crucially on local compactness,
and 
the result is no longer true in infinite dimensions 
as observed by Grossman \cite{grossman} and Atkin \cite{atkinhopfrinow}. 
The former constructed an infinite-dimensional ellipsoid in the space $\ell^2$ of 
square-summable sequences with points which cannot be connected by 
a minimal geodesic in the induced metric from $\ell^2$,
and the latter modified this construction to get points that cannot be joined by any geodesic at all. 
Interestingly, Ekeland~\cite{ekeland} showed that on a complete Riemannian Hilbert manifold 
the set of points attainable from a given one with a minimizing geodesic 
contains a dense $G_\delta$ 
(i.e., a countable intersection of open sets). 

One situation in which everything works nicely is 
%$\overline{\Omega(M)}^{_{H^1}}=H^1(S^1, M)$, 
the $H^1$ completion of the space of smooth loops 
$\overline{\Omega M}^{_{H^1}}$. 
Unlike many other examples discussed here, 
it is a genuine (strong) Riemannian Hilbert manifold in the topology generated by 
the distance function of \eqref{H1loopspace},
and by the result of 
Eliasson \cite{eliasson} any two of its points can be joined by a minimizing geodesic. 
%and the fact that it involves control of at least one derivative implies enough compactness 
%to get global results. 
%(the result is not valid for the loop space in the $L^2$ metric, for example~\textcolor{red}{Verify?}).

%%%%%%%%%%%%
\subsection{Exponential map and extendability of geodesics} 
If the Cauchy problem for the geodesic equation on a (possibly weak) Riemannian Hilbert manifold 
$\mathcal{M}$ is locally well-posed, then the exponential map of $\mathcal{M}$ can be defined as 
in finite dimensions. 
Using the scaling properties of geodesics we set 
\begin{align} \label{expdef} 
\exp_p : U \subset T_p\mathcal{M} \to \mathcal{M}, 
\qquad 
\exp_p(v) = \gamma(1),
\end{align} 
where $\gamma(t)$ is the unique geodesic 
from $\gamma(0)=p$ with initial velocity $\dot{\gamma}(0)=v$ in some open neighbourhood 
$U$ of zero in the tangent space at $p$. 

In general local (in time) well-posedness refers to constructing a unique solution 
for a given initial data on a short time interval which depends at least continuously on the data.
However, from the point of view of differential geometry, it is desirable if the dependence 
is at least $C^1$ smooth. 
Indeed, in this case applying the inverse function theorem for Banach manifolds, it is possible 
to deduce that (as in finite dimensions) the exponential map is a local diffeomorphism; this 
implies in particular nondegeneracy of the geodesic distance as in Section \ref{degeneratedistance}. 
Furthermore, other geometric tools such as Jacobi fields and curvature can be introduced 
to study rigorously stability as the problem of geodesic deviation 
(we shall elaborate on this in Section \ref{stability}).

The exponential maps 
defined on suitable Sobolev completions in the examples discussed so far 
are either at least $C^1$ smooth or else are continuous (even differentiable) but not $C^1$. 
The former include 
\begin{itemize}
\item 
the $L^2$ metric on the volumorphism group $\Diffmu^s(M)$ whose geodesics 
correspond to the Euler equations of ideal hydrodynamics; see \cite{Ebin-Marsden}, 
\item 
the $H^1$ metric on $\Diffmu^s(M)$ corresponding to the Lagrangian-averaged Euler equation;
see \cite{shkoller}, 
\item 
the $H^1$ metric on $\Diff^s(S^1)$ corresponding to the Camassa-Holm equation in \cite{ck} 
and its generalization to $\Diff^s(M)$ and the EPDiff equation in \cite{MP}, 
\item 
the $H^1$ metric on the free loop space $\overline{\Omega M}^{_{H^1}}$ in \cite{misiolekfreeloop}, 
\item 
the homogeneous $\dot{H}^1$ metric on $\Diff^s(S^1)/S^1$ corresponding to 
the Hunter-Saxton equation in \cite{K-M} and its generalization \eqref{HSmultid} 
on $\Diff^s(M)/\Diffmu^s(M)$; see \cite{klmp},
\item 
the right-invariant $L^2$ metrics on $\Met^s(M)$, $\Metmu^s(M)$ and $\Vol^s(M)$ 
in \cite{ebinmetrics, fg}, 
\item 
the noninvariant $L^2$ metric on $H^s(N, M)$ or $\Diff^s(M)$ whose geodesics 
are described by pointwise geodesics on $M$; see \cite{Ebin-Marsden}. 
\end{itemize} 
The metrics for which $C^1$ dependence fails include 
\begin{itemize}
\item 
the $L^2$ metric on the Virasoro group whose geodesic equation corresponds to 
the Korteweg-de Vries equation; see \cite{ckkt}, 
\item 
the right-invariant $L^2$ metric on $\Diff^s(S^1)$ which yields 
the (right-invariant) Burgers equation~\cite{ck} or its higher-dimensional generalization 
and the template-matching equation~\cite{MMvanishingdiffeo}, 
\item 
the $L^2$ metric on unit-parametrized curves in the plane yielding 
the whip equation \eqref{whip} in \cite{prestonwhip}, 
or on the equivalence classes of curves under reparametrizations in \cite{MM}.
\end{itemize}
Later on we will describe examples where the exponential map fails to be 
$C^1$ as a result of accumulation of conjugate points at $t=0$ 
(as in \cite{ck}, \cite{ckkt}, \cite{MKdV}, and \cite{prestonwhip}). 

%
%Finally we are interested in global existence for geodesics (i.e., a globally-defined exponential map). In finite dimensions this follows from metric completeness by the Hopf-Rinow theorem; as mentioned above the Hopf-Rinow theorem is not generally true in infinite dimensions. Alternatively one can prove global existence in finite dimensions if $M$ is a Lie group and the metric is right-invariant: this follows from the fact that the exponential map is defined for velocities smaller than some $\varepsilon$ at every point, and by right-invariance this $\varepsilon$ is the \emph{same} for every point. This would work in infinite dimensions if the metric were \emph{strong}, since then the Riemannian norm of velocity is conserved. However if our infinite-dimensional manifold has a weak metric, the norm of a tangent vector is conserved in the weak metric but may grow in the strong metric. Hence global existence requires an a priori bound on the strong norm of the velocity. For example, in ideal fluid mechanics on a manifold $M$, the $L^2$ norm of velocity is conserved, but we need a bound on the $H^s$ norm for $s>\frac{\dim{M}}{2} + 1$. It is sufficient to get a bound on the $C^1$ norm, and in fact it is sufficient to get a bound on the $L^{\infty}$ norm of vorticity, but whether this is possible depends on the detailed structure of the equations. 
%

It is well-known that the group of volumorphisms $\Diffmu^s(M)$ of a two-dimensional manifold $M$ 
equipped with the right-invariant $L^2$ metric is geodesically complete, 
i.e., its geodesics which correspond to solutions of the 2D incompressible Euler equations 
are defined globally in time when $s>2$.\footnote{We do not discuss weak solutions 
(for which much of the geometry seems to break down) of the PDE mentioned above.} 
This result is due to Wolibner \cite{w} with subsequent contributions by Yudovich \cite{y} 
and Kato \cite{kato} and follows from conservation of vorticity, although the argument is not routine. 
In three dimensions the problem is open and challenging. 

One might expect the right-invariant $H^1$ metric on the volumorphism group 
to be somewhat better behaved but as of now we have the same result: 
global existence in two dimensions is relatively easy~\cite{shkoller} but in three dimensions 
is unknown~\cite{houli}. 
On the other hand, it is known that smooth solutions of 
the one-dimensional Camassa-Holm equation \eqref{CH} break down for certain initial data 
\cite{mckeanblowup, mckeanlagrangian}. 
All solutions of the periodic Hunter-Saxton equation \eqref{HS1d} as well as its higher-dimensional 
generalization \eqref{HSmultid} are also known to blow up in finite time. 

For the non-invariant $L^2$ metric on $\Diff^s(M)$ (whose geodesics are given by 
pointwise geodesics on $M$) global existence clearly fails: 
for a typical initial velocity field two geodesics will eventually cross 
(which corresponds to a ``shock''). 
On the other hand in $H^s(M, M)$ geodesics exists for all time since such maps need not be injective. 

%Here part of the issue is the fact that solutions of the Burgers equation $u_t + \nabla_uu = 0$ blow up when the geodesics cross (this is known as a shock~\cite{KMshocks}), even though the flow is well-defined for all time; hence we have Eulerian blowup but global Lagrangian well-posedness, a consequence of the lack of smoothness of the group operation (which is used to get the Eulerian velocity of a geodesic tangent vector). The same sort of thing happens for the Camassa-Holm equation \eqref{CH}: we have global existence of geodesics in the $H^1$ metric on $C^{\infty}(S^1, S^1)$, although possibly not in $\Diff(S^1)$~\cite{mckeanlagrangian}; however the Camassa-Holm equation itself (the Eulerian equation of the geodesic flow) blows up in finite time if and only if the initial ``momentum'' $u-u_{xx}$ takes on both signs on $S^1$~\cite{mckeanblowup}. The Hunter-Saxton equation \eqref{HS1d} is known to blow up in finite time for any initial condition. On the other hand,  the geometric picture of it as a subset of the round sphere gives a way of defining global conservative solutions in $C^{\infty}(S^1, S^1)$ where the solution may leave the diffeomorphism group but it remains a homeomorphism~\cite{lenellsglobal}. The same is true for the multidimensional generalization \eqref{HSmultid} on the homogeneous space $\Diff(M)/\Diffmu(M)$, although it is not known whether there is actually a global lift to $C^{\infty}(M, M)$~\cite{klmp}.

Geodesics in the space of metrics $\Met^s(M)$ and volume forms $\Vol^s(M)$ 
typically become degenerate in finite time, see \cite{clarke} and \cite{fg}, 
while geodesics in $\Metmu^s(M)$ persist for all time \cite{ebinmetrics}. 
Both the $L^2$ and $H^1$ metrics on the homogeneous space of equivalence classes of curves 
in the plane admit geodesics that degenerate to points \cite{MM}. 
On the space of unit-speed curves, the whip equation would be expected to blow up in finite time 
physically~\cite{prestonwhip}, but this is not yet proved.

%%%%%%%%%%%%%%%%%%%%%%%%%%%%%
%%%%%%%%%%%%%%%%%%%%%%%%%%%%%%%%%%%%%%%%%%

\section{Jacobi fields, curvature, and stability} 
\label{stability} 
\nequation 

As mentioned in the Introduction, one application of Riemannian techniques 
in the study of equations of fluid dynamics has been to the problem of (Lagrangian) stability 
using the equation of geodesic deviation (the Jacobi equation)
which involves the curvature tensor. 
In this section we describe this approach for a general infinite-dimensional manifold 
equipped with a possibly weak Riemannian metric but whose exponential map is assumed 
to be at least $C^1$. 
We will discuss Jacobi fields (as infinitesimal perturbations) and the role played by 
sectional curvature and its sign. 
Various results for specific examples mentioned in the previous sections 
will be the subject of Section \ref{subsec:knownresults}. 

Let $\mathcal{M}$ be a (possibly weak) Riemannian Hilbert manifold 
whose geodesic equation is written in the form 
$$ 
\frac{D}{dt} \frac{d\gamma}{dt} = 0, 
$$
where $\frac{d\gamma}{dt}$ is the tangent vector field and $\frac{D}{dt}$ is the covariant derivative 
along the curve $\gamma(t)$ in $\mathcal{M}$. 
If $\bar{\gamma}(s,t)$ is a family of geodesics with $\bar{\gamma}(0,t)=\gamma(t)$ 
then the formula 
\begin{equation} \label{jacobidef} 
J(t) = \frac{\partial \bar{\gamma}}{\partial s}(0,t) 
\end{equation} 
gives a Jacobi field $J(t)$ along $\gamma$, i.e. a solution of the Jacobi equation  
\begin{equation} \label{jacobieqn} 
\frac{D^2J}{dt^2} + R\Big( J, \frac{d\gamma}{dt} \Big) \frac{d\gamma}{dt} = 0 
\end{equation} 
obtained by differentiating the geodesic equation in $s$ and evaluating at $s=0$. 
As in finite dimensions the Riemann curvature tensor $R$ of $\mathcal{M}$ 
arises here due to the fact that covariant derivatives do not commute in general. 
Furthermore, the basic result of Cartan applies as well so that 
for any $v$ and $w \in T_p\mathcal{M}$ we have 
\begin{equation} \label{expdiff} 
(D\exp_p)_{tv}(tw) = J(t) 
\end{equation} 
where $J(t)$ is the Jacobi field along $\gamma(t)=\exp_p(tv)$ solving \eqref{jacobieqn} 
with initial conditions $J(0)=0$ and $J'(0)=w$. 

Recall that the sectional curvature in the direction of the 2-plane spanned by the vectors 
$X_p$ and $Y_p \in T_p\mathcal{M}$ is given by the formula 
\begin{equation}\label{sectionalcurvature}
K(p) 
= 
\frac{\langle R(X_p, Y_p) Y_p, X_p\rangle}{\lvert X_p\rvert^2 \lvert Y_p\rvert^2 -\langle X_p, Y_p\rangle}.
\end{equation}

\begin{example} \upshape
On a two-dimensional Riemannian manifold the Jacobi equation can be reduced to 
a single ODE for a function $j(t)$ representing the component of $J$ orthogonal to 
$\dot{\gamma}$, which takes the form 
\begin{equation} \label{2djacobi}
\frac{d^2j(t)}{dt^2} + K(\gamma(t)) j(t) = 0,
\end{equation}
where $K$ is the sectional curvature at a point $\gamma(t)$. 
In the special case where $K$ is constant the solution of \eqref{2djacobi} with $j(0)=0$ is 
$$ 
j(t) = j'(0) \cdot \begin{cases} 
\frac{1}{\sqrt{K}} \sin{\sqrt{K}t} & K>0 \\
t & K=0 \\
\frac{1}{\sqrt{\lvert K\rvert}} \sinh{\sqrt{\lvert K\rvert} t} & K<0
\end{cases}.
$$
This simple special case is the source of much of our intuition about curvature and stability. 
Suppose that we know precisely the initial position of a particle traveling along a geodesic 
and its initial velocity only approximately. If $K>0$ then all Jacobi fields are bounded 
uniformly in time, and thus geodesics starting with nearby initial velocities 
will remain nearby for all time. 
If $K<0$ then the Jacobi fields grow exponentially in time, so that small errors are magnified 
and the motion is unpredictable. 
If $K=0$ then the growth is polynomial. 
On higher dimensional manifolds with variable curvature things become more subtle.
\end{example}

Recall that singular values of the Riemannian exponential map are called conjugate points. 
More precisely, two points $p$ and $q$ along a geodesic in $\mathcal{M}$ are conjugate 
if $D\exp_p$ viewed as a linear operator from $T_p\mathcal{M}$ to $T_q\mathcal{M}$ 
given by \eqref{expdiff} 
either fails to be injective (in which case the points are called \textit{mono-conjugate}) 
or it fails to be surjective (the points are called \textit{epi-conjugate}); see Grossman~\cite{grossman}. 
In finite dimensions both types coincide. 

Next, we state the Rauch comparison theorem for weak Riemannian metrics 
following Biliotti~\cite{biliotti}. 
This result relates growth of Jacobi fields to bounds on the sectional curvature 
and is a far-reaching generalization of Sturm's comparison theorem on oscillation of 
solutions of second order ODE. 

\begin{theorem} \label{rauchcomparison}
Let $\mathcal{M}$ and $\tilde{\mathcal{M}}$ be (possibly infinite-dimensional) 
weak Riemannian manifolds modeled on Hilbert spaces $E$ and $\tilde{E}$, 
with $E$ isometric to a closed subspace of $\tilde{E}$.
Assume that $\mathcal{M}$ and $\tilde{\mathcal{M}}$ have smooth Levi-Civita connections 
(and hence smooth exponential maps) 
with sectional curvatures $K$ and $\tilde{K}$. 
Let $\gamma$ and $\tilde{\gamma}$ be two geodesics of equal length 
and suppose that for every $X\in T_{\gamma(t)}\mathcal{M}$ 
and 
$\tilde{X} \in T_{\tilde{\gamma}(t)}\tilde{\mathcal{M}}$ 
$$ 
K(X, \gamma'(t)) \le \tilde{K}(\tilde{X}, \tilde{\gamma}'(t))\,.
$$
Let $J$ and $\tilde{J}$ be the Jacobi fields along $\gamma$ and $\tilde{\gamma}$ such that 
\begin{itemize}
\item $J(0)=0$ and $\tilde{J}(0)=0$,
\item $J'(0)$ is orthogonal to $\gamma'(0)$ and $\tilde{J}'(0)$ is orthogonal to $\tilde{\gamma}'(0)$, 
and
\item $\lVert J'(0)\rVert = \lVert \tilde{J}'(0)\rVert$.
\end{itemize}
If $\tilde{J}(t)$ is nowhere zero in the interval $(0, a]$ and if $\tilde{\gamma}$ has at most 
a finite number of points which are epi-conjugate but not mono-conjugate in $(0,a]$, 
then 
\begin{equation}\label{jacobicomparison}
\lVert J(t)\rVert \ge \lVert \tilde{J}(t)\rVert \quad \text{for all $t\in [0,a]$}.
\end{equation}
\end{theorem}

It often happens that such pathological points which are epi-conjugate but not mono-conjugate 
can fill out a whole interval 
(this is the case for the volumorphism group of a three-dimensional manifold~\cite{prestonWKB}) 
so that the criterion above may only be useful if the exponential map is \emph{Fredholm}; 
see Remark \ref{rem:expFred} below. 

\begin{remark} \upshape 
Theorem \ref{rauchcomparison} implies that if $\tilde{K}(X,\dot{\gamma})\ge -k$ 
for some positive constant $k$ 
(i.e., take $M$ to be a constant negative curvature space) 
and if $\tilde{\gamma}$ is free of conjugate points, then any Jacobi field along $\tilde{\gamma}$ 
in $\tilde{M}$ satisfies 
$\lVert \tilde{J}(t)\rVert \le \lVert J'(0)\rVert k^{-1/2} \sinh{kt}$, 
which gives essentially the maximum Lyapunov exponent for the system. 
In the opposite direction, if $K(X, \dot{\gamma})\le 0$ 
(i.e., take $\tilde{M}$ to be a flat space) then $\lVert J(t)\rVert \ge \lVert J'(0)\rVert t$, 
which can be interpreted as a weak instability with perturbations growing at least linearly in time.
In general, however, one should be cautious when drawing conclusions 
based on the Rauch theorem and one's finite-dimensional intuition: 
positive curvature does not imply stability, while negative curvature does not necessarily imply exponential instability, as we discuss below. 
\end{remark} 

For a general Riemannian manifold without additional structure one does not expect 
more precise results on the relation between curvature and stability. 
However, most of our examples have a group structure under which the Riemannian metric 
is right-invariant, which can be used to get additional information. 
To this end it will be convenient to decouple the Jacobi equation \eqref{jacobieqn} 
into two first-order equations. 

Namely, let $\mathcal{M}$ be a group $G$ with a right-invariant (weak) Riemannian metric. 
As in Section \ref{groupsubsection} its geodesics $\gamma(t)$ can be described by a pair of equations 
consisting of the Euler-Arnold equation \eqref{utBuu} and the flow equation \eqref{floweq} 
\begin{equation} \label{generalsplit} 
\frac{d\gamma}{dt} = DR_{\gamma}(u), \qquad \frac{du}{dt} + \ad_u^{\ast}u = 0 
\end{equation} 
defined in $G$ and $T_eG$ respectively. 
Let $\bar{\gamma}(s,t)$ be a family of such geodesics with $\bar{\gamma}(0,t) =\gamma(t)$ 
and with Eulerian velocity 
$\bar{u}=DR_{\bar{\gamma}^{-1}}\frac{d\bar{\gamma}}{dt}$. 
Setting $y = DR_{\gamma^{-1}}J$ 
(where $J$ is the Jacobi field along $\gamma$ as in \eqref{jacobidef}) 
and $z=\frac{\partial\bar{u}}{\partial s}|_{s=0}$ 
and differentiating both equations in \eqref{generalsplit} with respect to $s$ at $s=0$ 
we obtain a splitting of the Jacobi equation \eqref{jacobieqn} into 
\begin{align}
&\frac{dy}{dt} - \ad_uy = z   \label{linearizedflow} \\
&\frac{dz}{dt} + \ad_u^{\ast}z + \ad_z^{\ast}u = 0.   \label{linearizedeuler}
\end{align}

The \emph{linearized Euler equation} \eqref{linearizedeuler} can be used 
to define a notion of stability: 
a solution $u$ of the Euler equation is \emph{(linearly) stable} if every perturbation $z$ 
is bounded uniformly in time.
Using \eqref{linearizedflow} one can then relate this notion to 
stability of Lagrangian trajectories (as in \cite{P1} for the volumorphism group) 
and draw sharper conclusions about the behaviour of geodesics than is generally possible 
using only Rauch's theorem. The next two examples illustrate the subtleties. 

\begin{example}[\textit{Rigid body motion}] \upshape
Let $G=SO(3)$. Its Lie algebra $\mathfrak{so}(3)=T_eG$ 
is spanned by the vectors $e_1$, $e_2$ and $e_3$ satisfying 
$[e_1, e_2] = e_3$, $[e_2, e_3]=e_1$, $[e_3, e_1]=e_2$. 
This is the group of antisymmetric matrices represented as 
$$
x^1e_1 + x^2 e_2 + x^3 e_3 = \left( \begin{matrix} 
0 & -x^1 & -x^2 \\ 
x^1 & 0 & -x^3 \\
x^2 & x^3 & 0\end{matrix}\right). 
$$
The shape of the rigid body determines a left-invariant\footnote{All of our equations so far 
which have been stated for right-invariant metrics apply to left-invariant metrics after 
a possible change of the sign.} 
Riemannian metric by the conditions 
$\langle e_1, e_1\rangle = \lambda_1$, $\langle e_2, e_2\rangle = \lambda_2$, 
and $\langle e_3, e_3\rangle = \lambda_3$ 
for some positive numbers $\lambda_1, \lambda_2, \lambda_3$. 
The left-invariant analogue of \eqref{generalsplit} reads as follows
$$ 
\frac{d\gamma}{dt} = \gamma u, \qquad \frac{du}{dt} = \ad_u^{\ast}u 
$$
with the Euler-Arnold equation given explicitly by  
$$ 
\frac{du^1}{dt} = \frac{\lambda_2-\lambda_3}{\lambda_1} \, u^2 u^3, 
\qquad 
\frac{du^2}{dt} = \frac{\lambda_3-\lambda_1}{\lambda_2} \, u^1 u^3, 
\qquad 
\frac{du^3}{dt} = \frac{\lambda_1-\lambda_2}{\lambda_3} \, u^1 u^2.
$$
Consider one steady solution of this equation given by $u^1=u^3=0$ with $u^2=1$, 
supposing that $0<\lambda_1<\lambda_2<\lambda_3$. 
The linearized Euler equation \eqref{linearizedeuler} takes the form 
$$ 
\frac{dz^1}{dt} = \frac{\lambda_2-\lambda_3}{\lambda_1} \, z^3, 
\qquad 
\frac{dz^2}{dt} = 0, 
\qquad 
\frac{dz^3}{dt} = \frac{\lambda_1-\lambda_2}{\lambda_3} \, z^1 
$$
and if $\lambda_1<\lambda_2<\lambda_3$ it admits exponentially growing 
solutions.\footnote{This corresponds to the well-known fact that rotation of a rigid body 
about its largest and smallest axes is stable but rotation about the middle axis is unstable.}
The linearized flow equation \eqref{linearizedflow} takes the form 
$$ 
\frac{dy^1}{dt} + y^3 = z^1, 
\qquad 
\frac{dy^2}{dt} = z^2, 
\qquad 
\frac{dy^3}{dt} - y^1 = z^3, 
$$
so that $y(t)$ grows exponentially if $z(t)$ does.
On the other hand, 
sectional curvature in all directions containing $e_2$ can be made positive 
by a suitable choice of $\lambda$'s. Using e.g., the formulas of Milnor~\cite{milnor} 
we compute for $x=x^1e_1 + x^2 e_2$ that 
\begin{multline*} 
\langle R(e_2, x)x, e_2\rangle = \frac{(\lambda_2-\lambda_1)^2-\lambda_3^2 + 2\lambda_3(\lambda_2-\lambda_3+\lambda_1)}{4\lambda_3} \, (x^1)^2
\\+ \frac{(\lambda_3-\lambda_2)^2-\lambda_1^2 + 2\lambda_1(\lambda_1-\lambda_3-\lambda_2)}{4\lambda_1} \, (x^3)^2,
\end{multline*}
which is positive-definite if e.g.,
$\lambda_1=\frac{4}{5}$, $\lambda_2=1$, $\lambda_3=\frac{6}{5}$. 
Hence, we have positive curvature along the geodesic but exponentially growing Jacobi fields. 
This happens because the Rauch comparison theorem bounds Jacobi fields 
only up to the first conjugate point; beyond that point all bets are off.
\end{example}

\begin{example}[\textit{Couette flow}] \upshape 
Consider $\mathcal{M}=\Diffmu([0,1]\times S^1)$ with the $L^2$ metric. 
The ``plane-parallel Couette flow''  $u(x,y) = x \, \frac{\partial}{\partial y}$ 
is a steady solution of the Euler equation. 
This solution is known to be Eulerian stable~\cite{arnold72, orr} 
even though the sectional curvature is non-positive in all sections and typically negative~\cite{MStab}. 
A closer inspection reveals that the growth of all Jacobi fields is precisely linear~\cite{P1}. 
Hence, we do no better than what Rauch's theorem says: 
negative curvature need not imply exponential instability. 
\end{example} 

It is worth pointing out that one can relate Eulerian stability to Lagrangian stretching, 
i.e., to the growth of $\Ad_{\eta(t)}$ in the operator norm. 
Using the formula 
$d(\Ad_{\eta}y)/dt = \Ad_{\eta}( dy/dt - \ad_uy)$ 
and defining $Y$ and $Z$ by $y=\Ad_{\eta}Y$ and $z=\Ad_{\eta}Z$,
equation \eqref{linearizedflow} can be rewritten as 
\begin{equation} \label{adstarlinflow}
\frac{dY}{dt} = Z,
\end{equation} 
while equation \eqref{linearizedeuler} becomes 
\begin{equation} \label{adstarlineuler}
\frac{d}{dt} (\Ad_{\eta}^{\ast}\Ad_{\eta} Z) + \ad_Z^{\ast}u_0 = 0
\end{equation}
after incorporating conservation of vorticity $\Ad_{\eta}^{\ast}u = u_0$; see \cite{MP}. 
Observe that the operator $Z \to \Ad_{\eta(t)}^{\ast}\Ad_{\eta(t)}Z$ is selfadjoint and positive-definite 
while $Z \mapsto \ad_Z^\ast u_0$ is anti-selfadjoint with constant coefficients. 
This makes  \eqref{adstarlineuler} somewhat simpler to analyze than \eqref{linearizedeuler} 
(even if $u$ is independent of time) because $z \mapsto \ad_u^\ast z + \ad_z^\ast u$ 
is not selfadjoint; see \cite{vishikspectrum}. 

\begin{remark}[\textit{Fredholm exponential maps and conjugate points}] 
\label{rem:expFred} 
\upshape 
If the exponential map of a weak Riemannian manifold $\mathcal{M}$ is known to be smooth, 
one can ask about the distribution and nature of its singular values (conjugate points): 
can mono-conjugate and epi-conjugate points coincide, have finite multiplicity, 
or be discretely distributed along finite geodesic segments? 

These questions turn out to have positive answers if the derivative of the exponential map 
is a Fredholm operator between the tangent spaces of $\mathcal{M}$ 
with index zero.\footnote{In this case the exponential map is said to be a nonlinear Fredholm map 
of index zero (provided that $\mathcal{M}$ is connected).} 
This was first established for the free loop space $\overline{\Omega M}^{_{H^1}}$ 
with the Sobolev $H^1$ metric \eqref{H1loopspace} in \cite{misiolekfreeloop}. 
In this case the proof of Fredholmness is based on the fact that the curvature operator $R$ 
in \eqref{jacobieqn} is compact. 
In general one does not expect compactness. However, in the special case 
when $\mathcal{M}$ is a group $G$ with a right-invariant metric 
one can analyze the derivative of the exponential map using the pair of equations 
\eqref{adstarlinflow}--\eqref{adstarlineuler} 
to conclude that it is a sum of two terms, 
the first determined by the positive-definite operator $Z \to \Ad_\eta^\ast\Ad_\eta Z$ 
and the second a composition of a bounded map with $Y \to \ad_Y^\ast u_0$. 
The former is invertible. Thus, if the latter is compact then the exponential map 
will be Fredholm of index zero. 
This strategy works on the volumorphism group $\Diffmu^s(M)$ with the $L^2$ metric 
if $n=2$ (but not if $n \geq 3$) 
as well as for right-invariant Sobolev metrics of sufficiently high order 
on any diffeomorphism group; see \cite{MP}.

Fredholmness fails for the exponential map of the ellipsoid in $\ell^2$ with the induced metric 
(there are sequences of mono-conjugate points accumulating at an epi-conjugate point, 
or a mono-conjugate point of infinite order~\cite{grossman}). 
It also fails for the free loop space with the $L^2$ metric~\cite{misiolekfreeloop} 
and for the volumorphism group $\Diffmu^s(M)$ of a three-dimensional manifold 
in the $L^2$ metric~\cite{prestonWKB}. 
In the latter case, we can have mono-conjugate points that are dense in an interval 
and epi-conjugate points that fill up an interval. 

On the other hand, 
smoothness of the exponential map implies (by the inverse function theorem) 
that any sufficiently short geodesic segment is free of conjugate points 
and hence is locally minimizing. 
In particular, 
if for some $t_n\searrow 0$ the points $\gamma(t_n)$ along a given geodesic are 
mono-conjugate to $\gamma(0)$ then the exponential is not $C^1$. 
This method was used to prove that the exponential maps associated to 
the KdV equation~\cite{MKdV}, the right-invariant Burgers equation~\cite{ck}, 
and the whip equation \eqref{whip}~\cite{prestonwhip} cannot have 
$C^1$ exponential maps. 
\end{remark}

\section{The sign of the curvature: previous results} 
\label{subsec:knownresults}
\nequation 

In the remainder of the paper we will focus on the sign of the sectional curvature 
in the examples described above. 
It turns out that with few exceptions sectional curvature can be positive or negative 
depending on the two-dimensional direction. 
This section contains a survey of known results and techniques. 
New results will be presented in Sections \ref{circlediffeo} and \ref{bothsignsnow}. 

%\begin{remark}\upshape
%It is quite difficult to work with a general vector field on an infinite-dimensional manifold, since e.g., on $\Diff(S^1)$ such a field would be described in terms of nonlinear functionals of functions on $S^1$. Thus usually the covariant derivative $\frac{D}{dt}$ along a particular curve is easier to work with (corresponding to a time-dependent vector field on $S^1$, which is much more concrete). But the fact that curvature depends only on the values of the vector fields at a point implies that one can do explicit calculations given a basis of vector fields. Such a basis is provided in concrete cases by the following simple construction. Recall that on the space of manifold maps $C^{\infty}(N, M)$ a tangent vector to $\eta$ is a map $u\colon N\to TM$ such that $u(p)\in T_{\eta(p)}M$ for every $p\in N$. Given any vector field $v$ on $M$, the map $v\circ\eta$ satisfies this condition, and hence the map $V\colon C^{\infty}(N,M)\to TC^{\infty}(N, M)$ given by $V_{\eta} = v\circ\eta$ is a vector field on $C^{\infty}(N, M)$. One can obtain all vectors at a given $\eta$ in this way if $\eta$ is an embedding of $N$ into $M$. In particular on $\Diff(M)$ such vector fields are right-invariant, and they span every tangent space. 
%\end{remark}

%%%%%
%\subsubsection{}
%%%
The simplest curvature formula arises on $\DiffM \subset C^{\infty}(M, M)$ 
equipped with the $L^2$ metric \eqref{noninvariantL2}. 
If $U=u\circ\eta$ and $V=v\circ\eta$ are two vector fields on $\DiffM$ where $u, v \in T_e\DiffM$ 
then the covariant derivative of the $L^2$ metric is computed in terms of 
the covariant derivative on $M$ as 
$(\nabla_UV)_{\eta} = (\nabla_uv)\circ\eta$. 
Therefore, the $L^2$ curvature of $\DiffM$ is completely determined by 
the Riemannian curvature of $M$ 
and can be computed directly from the definition 
$$ 
\llangle R(u\circ\eta, v\circ\eta)v\circ\eta, u\circ\eta\rrangle_{L^2} 
= 
\int_M \langle R(u,v)v,u\rangle\circ\eta\, d\mu. 
$$
A few simple observations can be made based on this formula. 
If the vector fields $u$ and $v$ have unit $L^2$ norms and 
are chosen to have disjoint supports then the integral on the right-hand side will be zero. 
Consequently, the $L^2$ curvature of $\DiffM$ cannot be strictly positive or strictly negative 
even if $M$ has constant (non-zero) curvature.  
Furthermore, it is also clear that it can be non-negative (non-positive) 
if and only if 
the sectional curvature of $M$ is non-negative (non-positive). 

One can similarly obtain relatively simple expressions for the $L^2$ curvature of 
the volumorphism group $\Diffmu(M)$ (following \cite{MStab}) 
or the free loop space $\Omega M$ 
as well as the curvature of the unit-speed loops in $\mathbb{R}^2$ (see \cite{prestonwhip}).
To do this we use submanifold geometry.

%%%%%%%%%
%\subsubsection{} 
%%%%
As in finite dimensional geometry, 
if $\mathcal{N}$ is a submanifold of a (weak) Riemannian manifold $\mathcal{M}$ 
then the induced Levi-Civita connection on $\mathcal{N}$ 
is related to that on $\mathcal{M}$ via the second fundamental form 
\begin{equation}\label{secondfundy}
\Pi(U,V) = \nabla^{\mathcal{M}}_UV-\nabla^{\mathcal{N}}_UV
\end{equation} 
where $U$ and $V$ are vector fields tangent to $\mathcal{N}$. 
$\Pi$ is symmetric and tensorial 
(i.e.,  its value at any $p\in \mathcal{N}$ depends only on the values $U_p$ and $V_p$) 
and 
the curvature of $\mathcal{N}$ can then be computed using the Gauss-Codazzi formula 
\begin{equation}\label{gausscodazzi}
\begin{split}
\llangle R^{\mathcal{N}}(U,V)V, U\rrangle 
&= \llangle R^{\mathcal{M}}(U,V)V,U\rrangle
\\ 
&\qquad\qquad+ \llangle \Pi(U,U), \Pi(V,V)\rrangle - \llangle \Pi(U,V), \Pi(U,V)\rrangle. 
\end{split}
\end{equation}

If $\mathcal{N}=\DiffmuM$ is the volumorphism group and $\mathcal{M}=\DiffM$ is the group 
of all diffeomorphisms with the $L^2$ metric 
then the corresponding second fundamental form is 
$\Pi(U,V) = \grad \Laplacian^{-1} \diver{(\nabla_uv)}\circ\eta$ 
where 
$U=u\circ\eta$, $V=v\circ\eta$ and where $u$ and $v$ are divergence free vector fields on $M$. 
The following theorem summarizes the known results in this important case;
see \cite{lukatsky90, lukatsky, MStab, PNonpositive, rouchon, sm} 
\begin{theorem} \label{thm:51} 
Let $M$ be a compact manifold of dimension $n \geq 2$, possibly with boundary. 
Consider the volumorphism group $\Diffmu(M)$ with the $L^2$ metric. 
For any $u\in T_e\Diffmu(M)$ define 
$$ 
K_{\min}(u) 
= 
\inf_{v\in T_e\Diffmu(M)} 
\frac{\llangle R(v,u)u,v\rrangle_{L^2}}{\lVert v\rVert^2_{L^2}\lVert u\rVert^2_{L^2} 
- 
\llangle u,v\rrangle^2_{L^2}} 
$$
to be 
the minimum sectional curvature in directions containing $u$, 
and similarly define $K_{\max}(u)$ to be the maximum curvature. 
Then we have
\begin{enumerate}
\item 
$K_{\min}(u) < 0$, 
unless $u$ is a Killing field in which case $K_{\min}(u)=0$; 
\item 
if $n \ge 3$ and $M$ is flat then $K_{\max}(u)>0$, 
unless $\diver{(\nabla_uu)}=0$ in which case $K_{\max}(u)=0$; 
\item 
if $n=2$ and $M$ is flat then $K_{\max}(u)>0$, 
unless $u$ is plane parallel $u = f(x) \, \frac{\partial}{\partial y}$ 
or purely rotational $u = f(r) \, \frac{\partial}{\partial \theta}$, 
in either of which cases $K_{\max}(u)=0$. 
\end{enumerate}
\end{theorem}
\begin{proof} 
All three statements follow from the Gauss-Codazzi formula (\ref{gausscodazzi}) 
which is the most effective way to determine the sign of the curvature of $\Diffmu(M)$. 
Part (1) can be found in \cite{rouchon} when $n=3$ and $M$ is flat although the technique works 
in any dimension. 
In part (2), the fact that $\diver{(\nabla_uu)}=0$ implies $K_{\max}(u)=0$ 
is essentially due to \cite{MStab} 
while the converse can be proved using the same approximation scheme as in \cite{Phardest} 
(used there to find conjugate points). 
Finally, part (3) is a special case of a result in \cite{PNonpositive} 
which works for any steady flow on a surface. 
\end{proof}

%%%%%%%
%\subsubsection{}
%%%%
In particular, it follows from Theorem \ref{thm:51} that for any $M$ of dimension $n \ge 2$ 
the $L^2$ sectional curvature of $\Diffmu(M)$ assumes both signs. 
The first examples in the special case of the flat 2-torus 
$M=\mathbb{T}^2$ were worked out by Arnold \cite{A}, 
who used the Lie-theoretic approach of Section \ref{groupsubsection}. 
He derived a general formula for the sectional curvature of 
a group $\mathcal{M}=G$ with a right-invariant metric in terms of 
the coadjoint operator \eqref{opB} as 
\begin{equation} \label{generalcurvature}
\begin{split}
\llangle R(u,v)v,u\rrangle 
= 
\frac{1}{4} \lVert \ad^\ast_v u &+ \ad^\ast_u v \rVert^2 
- 
\llangle \ad^\ast_u u, \ad^\ast_v v \rrangle 
\\ 
&- 
\frac{3}{4} \lVert \ad_uv \rVert^2 
+ 
\frac{1}{2} \llangle \ad_uv, \ad^\ast_v u - \ad^\ast_u v \rrangle 
\end{split}
\end{equation}
and applied it to the case of $\Diffmu(\mathbb{T}^2)$. 

\begin{example} \upshape 
If $u$ and $v$ are vector fields on $\mathbb{T}^2$ with the stream functions 
$f(x,y) = \cos{(jx+ky)}$ and $g(x,y) = \cos{(lx+my)}$ 
then the (unnormalized) sectional curvature of $\Diffmu(\mathbb{T}^2)$ is 
\begin{equation} \label{arnoldeq}
\llangle R(u,v)v,u\rrangle_{L^2} 
= 
-\frac{\pi^2(jm-kl)^4(j^2+k^2+l^2+m^2)}{\big((j+l)^2+(k+m)^2\big)\big((j-l)^2+(k-m)^2\big)} < 0. 
\end{equation}
On the other hand, if we pick 
$f(x,y) = \cos{(3kx-y)}+\cos{(3kx+2y)}$ 
and 
$g(x,y)=\cos{(kx+y)}+\cos{(kx-2y)}$ 
then 
$$ 
\lim_{k\to\infty} K(u,v) = \frac{9}{8\pi^2}. 
$$
Thus it is easy to find negative curvature, but there are many sections with positive curvature as well.
\end{example}

%%%%%%%%%%%%
%\subsubsection{}
%%%%%
There are two cases in which the curvature is known to have a remarkably simple form. 
The first is $\Diff(S^1)$ with the right-invariant $L^2$ metric whose curvature at the identity 
(and hence everywhere by right invariance) is given by 
\begin{equation}\label{1dburgerscurvature}
\llangle R(u,v)v,u\rrangle_{L^2} = \int_{S^1} ( uv_x - vu_x )^2 \, dx 
\end{equation}
and thus is non-negative. 
We will show in Section \ref{bothsignsnow} that this attractive formula does not generalize 
to higher dimensions and the sectional curvature of the right-invariant $L^2$ metric 
on $\DiffM$ can assume both signs.

The second case is 
the space of densities $\Diff(M)/\Diffmu(M)$ equipped with the homogeneous 
Sobolev $\dot{H}^1$ metric obtained by setting $b=1$ and $a=c=0$ in \eqref{rightinvmetric}. 
This space turns out to be isometric to the round sphere of radius $2$ 
and therefore has constant positive curvature 
\begin{equation} \label{H1K1/4} 
\llangle R(u,v)v, u \rrangle_{\dot{H}^1} 
= 
\frac{1}{4} \Big( \|u\|^2_{\dot{H}^1} \|v\|^2_{\dot{H}^1} - \llangle u, v \rrangle_{\dot{H}^1}^2 \Big). 
\end{equation} 
We refer to \cite{L, klmp} for detailed calculations. 

When the exponential map is not smooth, the curvature may be positive but unbounded above; this allows for conjugate points that occur arbitrarily close to any given point along a geodesic. 
In such situations 
one cannot determine stability studying geodesic deviation even for short times 
since the Rauch Theorem \ref{rauchcomparison} immediately fails. 
%Sometimes, as in the case of the whip equation \eqref{whip}, 
%it may be possible to use the splitting 
%\eqref{linearizedflow}--\eqref{linearizedeuler} 
%and obtain bounds on Jacobi fields. 
This applies for example to the right-invariant $L^2$ metric on the Virasoro group 
(which yields the KdV equation) 
whose exponential map is also known not to be smooth and whose 
sectional curvature is unbounded and of both signs~\cite{MKdV, ckkt}. 

\begin{remark} \upshape 
In what follows we will often use the letter $S$ as a shorthand notation for 
the (non-normalized) sectional curvature 
$
S(u,v) = \llangle R(u,v)v,u \rrangle
$
if the metric used is clear from the context. 
\end{remark}

\section{The sign of the curvature: the one dimensional case} 
\label{circlediffeo} 
\nequation 

In this section we present new results on the sign of the sectional curvature 
of the right-invariant $H^1$ metric on the group $\Diff(S^1)$. 
In this case the Sobolev $H^1$-metric (\ref{rightinvmetric}) reduces to the $a$-$b$ metric 
\begin{equation} \label{abmetric}
\llangle u,v\rrangle_{H^1} 
= 
\int_0^1 \big( au v + bu_x v_x \big) \, dx
\end{equation}
with $a>0$ and $b>0$. The corresponding Euler-Arnold equation \eqref{utBuu} reads 
\begin{equation} \label{1Dgeodesiceq} 
m_t = -3 a uu_x +  b(2 u_xu_{xx}  +uu_{xxx}), 
\qquad 
m = Au = au - bu_{xx}.
\end{equation}
For $a=b=1$ we get the periodic Camassa-Holm equation with period 1. 
For other values we can rescale by $y=x\sqrt{\frac{a}{b}}$ and $s=t\sqrt{\frac{a}{b}}$ 
so that \eqref{1Dgeodesiceq} becomes
$$
u_s - u_{syy} + 3 uu_y - 2u_yu_{yy} - u u_{yyy} = 0, 
$$
which is the Camassa-Holm equation with period $\sqrt{a/b}$.

Recall from Section \ref{subsec:knownresults} that the sectional curvature of $\Diff(S^1)$ 
equipped with the $a$-$b$ metric \eqref{abmetric} is already known in the ``end-point'' cases 
where either $a=0$ or $b=0$. 
In the former case, the sectional curvature is positive and constant (see \eqref{H1K1/4}), 
and 
in the latter it is non-negative (see \eqref{1dburgerscurvature}). 
We will next show that when both $a$ and $b$ are positive, the curvature 
of $\Diff(S^1)$ can assume both signs. 

\begin{lemma} \label{1Dcurvprop} 
The sectional curvature of $\Diff(S^1)$ endowed with the right-invariant $a$-$b$ metric \eqref{abmetric} 
where $a > 0$ and $b > 0$ is given by 
\begin{equation} \label{abcurvature}
S(u,v) 
= 
\llangle R(u,v)v, u \rrangle_{H^1} 
= 
\llangle \Gamma(u,v), \Gamma(u,v)\rrangle_{H^1} 
- 
\llangle \Gamma(u,u), \Gamma(v,v)\rrangle_{H^1}, 
\end{equation}
for any $u,v \in T_e \Diff(S^1)$, where $\Gamma$ is the Christoffel map defined by
\begin{equation} \label{1DGammadef}
\Gamma(u, v) 
= 
A^{-1}\partial_x\left(auv + \frac{b}{2}u_x v_x\right), 
\qquad 
A = a - b \partial_x^2.
\end{equation}
\end{lemma}
\begin{proof}
We have 
$\ad^\ast_v u = A^{-1}(2auv_x+avu_x-2bv_xu_{xx}-bvu_{xxx})$ 
and therefore an easily verified identity 
$$
\ad^\ast_v u + \ad^\ast_u v = \partial_x(uv) + 2\Gamma(u,v),
$$
which yields 
$$
\ad^\ast_u u = \tfrac{1}{2} \partial_x(u^2) + \Gamma(u,u). 
$$ 

Using the general curvature formula \eqref{generalcurvature}, we obtain 
$$
\llangle R(u,v)v,u \rrangle_{H^1} 
= 
\lVert \Gamma(u,v)\rVert^2_{H^1} 
- 
\llangle \Gamma(u,u), \Gamma(v,v)\rrangle_{H^1} + L(u,v),
$$
where 
\begin{align*}
L(u,v) = \llangle \Gamma(u,v), \partial_x(uv)\rrangle_{H^1} 
&- 
\tfrac{1}{2} \llangle \Gamma(u,u), \partial_x(v^2)\rrangle_{H^1} 
- 
\tfrac{1}{2} \llangle \Gamma(v,v), \partial_x(u^2) \rrangle_{H^1}          \\
&+ 
\tfrac{1}{4} \lVert \partial_x(uv)\rVert^2_{H^1} 
- 
\tfrac{1}{4} \llangle \partial_x(u^2), \partial_x(v^2)\rrangle_{H^1}            \\
&- 
\tfrac{3}{4} \lVert \ad_uv\rVert^2_{H^1} 
+ 
\tfrac{1}{2} \llangle \ad_uv, \ad^\ast_v u - \ad^\ast_u v \rrangle_{H^1} 
\end{align*}
and $\ad_uv = -v_xu + u_x v$. 
Since 
$\llangle f, A^{-1}g\rrangle_{H^1} = \int_{S^1} fg\, dx$ 
for any functions $f$ and $g$, we can perform all these computations without ever 
explicitly evaluating $A^{-1}$. 
A lengthy computation involving integration by parts shows that $L(u,v)$ is always zero.
\end{proof}

The following proposition shows that it is easy to find sections of positive $H^1$ curvature 
on $\Diff(S^1)$ with the $a$-$b$ metric. In fact, the curvature is strictly positive 
along all subspaces spanned by two trigonometric functions.

\begin{theorem} \label{1Dpositivecurvprop}
Consider $\Diff(S^1)$ endowed with the right-invariant $H^1$ metric given at the identity by 
\eqref{abmetric} with $a > 0$ and $b > 0$.
If $k$ and $l$ are strictly positive distinct integer multiples of $2\pi$, then 
\begin{equation} \label{SSSC} 
S(\cos{kx}, \cos{lx}) = S(\cos{kx}, \sin{lx})= S(\sin{kx}, \sin{lx}) = C(k,l) >0,
\end{equation}
where $S(u,v) = \llangle R(u,v)v,u \rrangle_{H^1}$ and 
$$
C(k,l) 
= 
\frac{1}{8}\left(\frac{(a + \frac{b}{2}kl)^2}{a + b (k-l)^2} (k-l)^2
+ 
\frac{(a - \frac{b}{2}kl)^2}{a + b (k+l)^2}(k+l)^2 \right). 
$$
Moreover, for $k$ an integer multiple of $2\pi$, we have
\begin{align*} 
&S(\cos{kx}, \sin{kx}) = 2C(k,k) = \frac{(a - \frac{b}{2}k^2)^2}{a + 4b k^2} k^2 > 0,
	\\
&S(\cos{kx}, 1) = S(\sin{kx}, 1) = 2C(k,0) = \frac{a^2k^2}{2(a + b k^2)} >0.
\end{align*}	
\end{theorem}
\begin{proof}
Let $u=\cos{kx}$ and $v=\cos{lx}$. 
From the definition \eqref{1DGammadef} of the Christoffel map $\Gamma$ 
we have 
\begin{align*}
&\Gamma(u,u) 
= 
- \frac{2ak-bk^3}{2a+8bk^2} \sin{(2kx)}, 
\qquad\quad 
\Gamma(v,v) = -\frac{2al-bl^3}{2a+8bl^2} \sin{(2lx)},
	 \\
&\Gamma(u,v) 
= 
- \frac{(k+l)(2a-bkl)}{4a+4b(k+l)^2} \sin{(k+l)x} 
- 
\frac{(k-l)(2a+bkl)}{4a+4b(k-l)^2} \sin{(k-l)x}.
\end{align*} 
Substituting these into \eqref{abcurvature} gives the formula for $S(\cos{kx}, \cos{lx})$.
The other formulas are proved in a similar way.
\end{proof}

Sections of negative $H^1$ curvature are trickier to find. 

\begin{theorem} \label{abnonzero1d}
For every choice of $a >0$ and $b> 0$, 
there exist velocity fields $u$ and $v$ such that the sectional curvature of $\Diff(S^1)$ 
endowed with the $a$-$b$ metric \eqref{abmetric} is strictly negative, i.e.,
$S(u,v)=\llangle R(u,v)v,u \rrangle_{H^1} <0$.
\end{theorem} 
\begin{proof} 
Set $\alpha = a/(4\pi^2b)$ and pick 
$$ 
u(x) = \phi + \cos{4\pi x}, 
\qquad
v(x) = \sin{2\pi x},
$$
where
$$
\phi =-\frac{3}{2}\frac{\alpha^2-\alpha-2}{\alpha(\alpha+4)};
$$
substitution into \eqref{abcurvature} yields
$$ 
S(u,v) 
= 
\frac{2b\pi^4 (\alpha^4+18\alpha^3 
+ 
357\alpha^2-20\alpha-36)}{(\alpha+9)(\alpha+4)^2}. 
$$
Observe that this quantity is negative for $0<\alpha \le 0.34$.

The second example is constructed differently and works when $\alpha \ge 0.34$. 
Choose a positive integer $j$ such that 
$\frac{1}{2}\sqrt{\alpha/0.34} < j \le \sqrt{\alpha/0.34}$ 
and define $r=\alpha/j^2$ so that $0.34 \le r < 1.36$. 
Set 
$$
\psi = \sqrt{-\frac{(73r^2-188r+45)(r+16)}{128(r+9)(r-2)^2}}. 
$$ 
It is easy to see that this is defined in the range specified above. 
Set 
$$
u(x) = \cos{2\pi jx} + \psi\cos{4\pi jx} 
\quad \mathrm{and} \quad 
v(x) = \sin{2\pi jx} + 2\psi\sin{4\pi jx}.
$$
Substituting into \eqref{abcurvature}, we obtain
$$
S(u,v) = -\frac{3\pi^4bj^4}{64} \frac{P(r)}{(r+9)^2 (r+4)(r-2)^2},
$$
where 
$P(r)
=
1435 r^6 + 21940 r^5 - 55074 r^4 - 222512 r^3 + 584323 r^2 - 215364 r + 15552. 
$
This quantity  $S(u,v)$ is negative for $0.34\le r<1.36$, as desired.
\end{proof}

\begin{remark}[\it The $\mu$CH equation] \label{mu-CH} \upshape
An interesting example of a right-invariant $H^1$-type metric on $\Diff(S^1)$ 
is given at the identity by 
$$ 
\llangle u, v\rrangle_{H^1_\mu} = c\, \mu(u) \mu(v) + \int_{S^1} u'(x)v'(x) \, dx\,,
$$
for any positive constant $c$, where $\mu(u):=\int_{S^1}u(x)dx$ 
is the mean value of the field over the circle. 
This metric yields yet another integrable evolution equation 
$$
u_{txx} - 2c\mu(u)u_x + 2u_x u_{xx} + uu_{xxx} = 0
$$
as a geodesic equation on the diffeomorphism group which ``interpolates'' between 
the Hunter-Saxton and Camassa-Holm equations.
(This equation is sometimes called the $\mu$HS or $\mu$CH equation.) 
The group $\Diff(S^1)$ equipped with the $H_\mu^1$ metric above admits sections of 
negative curvature; e.g., $S(u,v) < 0$ whenever 
$$ 
u(x) = \frac{3\pi^2 k^2}{c} + \cos{(4\pi k x)}, \qquad v(x) = \sin{(2\pi k x)}\,,
$$
and $k$ is any nonzero integer.
\end{remark}

%%%%%%%%%%%%%%%%%%%%%%%%%%%%%%%%%%%%%%%%%%%%
%%%%%%%%%%%%%%%%%%%%%%%%%%%%%%%%%%%%%%%%%%%%%%%%%

\section{The sign of the curvature: higher dimensions} 
\label{bothsignsnow} 
\nequation 

In order to simplify the formulas we present the results for the case when $M$ 
is the flat torus $\mathbb{T}^n=\mathbb{R}^n/\mathbb{Z}^n$. 
We recall the general formula for the Sobolev metric \eqref{rightinvmetric} 
in the form 
\begin{equation} \label{abcmetric}
\llangle u, v \rrangle_{H^1} 
= 
a \int_{\mathbb{T}^n} \langle u, v\rangle \, d\mu  
+  
b \int_{\mathbb{T}^n} \delta u^{\flat} \cdot \delta v^{\flat} \, d\mu
+ 
c \int_{\mathbb{T}^n} \langle du^\flat, dv^\flat\rangle \, d\mu 
\end{equation} 
and observe from \eqref{H1K1/4} that when $a=c=0$ then the corresponding 
sectional curvature of $\Diff(\mathbb{T}^n)$ is strictly positive and constant. 

In this section, we show that in the general $a$-$b$-$c$ case the sectional curvature of 
\eqref{abcmetric} on $\Diff(\mathbb{T}^n)$ assumes both signs. 
The case when at least two of the parameters $a, b, c$ are nonzero is treated in 
Section \ref{subsec:epdiff}, 
the case $b=c=0$ is treated in Section \ref{subsec:burgers}, 
and 
the case $a=b=0$ on the subgroup $\Diff_{\mu,\text{ex}}(\mathbb{T}^2)$ 
is in Section \ref{subsec:abzero}.
All 2D examples discussed below generalize naturally to higher dimensions.

%%%%%%%%%%%%%%%%%%%%%%%%%%%%%%

\subsection{The $H^1$-metric  on $\Diff(\mathbb{T}^n)$:  the EPDiff equation}
\label{subsec:epdiff}

In the case when all parameters $a$, $b$ and $c$ of the $H^1$ metric \eqref{abcmetric} 
are strictly positive the Euler-Arnold equation (\ref{utBuu}) 
is a multidimensional generalization of the Camassa-Holm equation. 
In the special case where $a=b=c=1$ and the manifold is 
a flat torus\footnote{In general, EPDiff involves 
the rough Laplacian $\nabla^*\nabla$ rather than 
the Hodge Laplacian $d\delta+\delta d$; these operators differ by 
a Ricci curvature term due to the Bochner-Weitzenb\"ock formula. 
If the manifold is Einstein the EPDiff metric is a special case of \eqref{abcmetric}.} 
we obtain the EPDiff equation \cite{HMR}.

\begin{theorem}\label{twononzero}
If $M$ is the flat torus $\mathbb{T}^n$ and at least two of the parameters $a,b,c$ 
are nonzero in \eqref{abcmetric}, then the curvature of $\Diff(\mathbb{T}^n)$ 
takes on both signs. 
\end{theorem}
\begin{proof}
The formula for the coadjoint operator for the $a$-$b$-$c$ metric (\ref{abcmetric}) is 
\begin{equation}\label{generalBop}
\ad^\ast_v u 
= 
A^{-1}\big( (\diver{v})Au + d\langle Au, v\rangle + \iota_vdAu\big),
\end{equation}
where $Av = av^\flat + b d \delta v^\flat + c \delta dv^b$; see \cite{klmp}.

If we pick 
$u = f(x) \tfrac{\partial}{\partial x}$ and $v = g(x) \tfrac{\partial}{\partial x}$ 
then \eqref{generalBop} gives 
$$ 
\ad^\ast_v u 
= 
(a - b\partial_x^2)^{-1} \Big( a\big( 2 g_x f + f_x g \big) - b\big( 2g_x f_{xx} + g f_{xxx} \big) \Big). 
$$
The value of $c$ is irrelevant in this case, since $du^{\flat} = dv^{\flat}=0$.
Thus, formula \eqref{1dburgerscurvature} yields examples with positive curvature 
whenever $a \neq 0$. 
If $a = 0$ but $b$ is nonzero, then positive curvature directions exist 
since this space is isometric to a sphere~\cite{L}. 
Finally, if $a$ and $b$ are both nonzero, then Proposition \ref{abnonzero1d} 
yields examples of negative curvature. 

To finish the proof we need negative-curvature examples 
when $b=0$ with $a$ and $c$ both nonzero, 
and 
when $a=0$ while $b$ and $c$ are both nonzero. 
We will present them for the two-dimensional flat torus $\mathbb{T}^2$.

Let $u = f(x)  \tfrac{\partial}{\partial y}$ and $v = g(x)\tfrac{\partial}{\partial x}$ 
so that 
$Au =(af - cf_{xx})\tfrac{\partial}{\partial y}$ 
and 
$Av =(ag - bg_{xx})\tfrac{\partial}{\partial x}$. 
If $a \ne 0$, then using \eqref{generalBop} we find 
%\begin{align*}
%A B(u,u) 
%&= \Big( -af(x) f'(x) + c f'(x) f''(x) \Big) \, \partial_x, \\
%A B(u,v) 
%&= \Big( -a[ g'(x) f(x) + g(x) f'(x)] 
%+ c[ g'(x) f''(x) + g(x) f'''(x)] \Big) \, \partial_y, \\
%A B(v,u) &= 0, \\
%A B(v,v) &= \Big( -3 a g(x) g'(x) 
%+ b [2g'(x)g''(x) + g(x) g'''(x)] \Big) \, \partial_x,
%\end{align*}
%where $A = (a + b d\delta + c \delta d)$. 
%If $a\ne 0$ we can write these more explicitly as 
%
\begin{align*}
\ad^\ast_u u 
&= 
(a-b\partial_x^2)^{-1}  
\Big( af f_x - c f_x f_{xx} \Big) \tfrac{\partial}{\partial x},
	\\
\ad^\ast_v u 
&= 
(a-c\partial_x^2)^{-1} \Big( a\big( g_x f + g f_x \big) 
- 
c\big( g_x f_{xx} + g f_{xxx} \big) \Big)  \tfrac{\partial}{\partial y},
	\\
\ad^\ast_u v 
&= 0,
	\\
\ad^\ast_v v 
&= 
(a-b\partial_x^2)^{-1}  
\Big( 3 a g g_x - b (2g_x g_{xx} + g g_{xxx} ) \Big)  \tfrac{\partial}{\partial x}.
\end{align*}
We also have $\ad_uv = g f_x \tfrac{\partial}{\partial y}$. 
If $f(x)=g(x)=\sin{kx}$ for some $k$ which is an integer multiple of $2\pi$, 
then it is easy to see that 
\begin{align*}
\ad^\ast_u u 
&= 
\frac{k(a+ck^2)}{2(a+4bk^2)}\, \sin{2kx} \, \tfrac{\partial}{\partial x}, \\
\ad^\ast_v u 
&= 
\frac{k(a+ck^2)}{a+4ck^2} \, \sin{2kx} \, \tfrac{\partial}{\partial y}, \\
\ad^\ast_v v 
&= 
\frac{3k(a+bk^2)}{2(a+4bk^2)} \, \sin{2kx} \, \tfrac{\partial}{\partial x}, \\
\ad_uv 
&= 
\frac{k}{2} \, \sin{2kx} \, \tfrac{\partial}{\partial y}.
\end{align*}

Note that these formulas are valid if $a=0$ as well, 
as long as $b\ne 0$ and $c\ne 0$. 
Indeed, in this case all vectors have to be projected to the orthogonal complement of 
the harmonic fields,\footnote{The metric \eqref{abcmetric} on $\Diff(\mathbb{T}^2)$ with $a=0$ is degenerate and only defined on the homogeneous space $\Diff(\mathbb{T}^2)/\mathbb{T}^2$. Hence everything is only defined modulo harmonic fields on $\mathbb{T}^2$.} but since 
$$
\int_0^1\int_0^1 \sin{kx} \, dx\,dy 
= 
\int_0^1\int_0^1 \sin{2kx} \, dx \, dy =0,
$$
we see that all components are already orthogonal to the harmonic fields.

Substitution into \eqref{generalcurvature} now gives the formula 
$$ 
S(u,v) 
= 
-\frac{k^2(7a^3-8a^2bk^2+56a^2ck^2+44ack^4b+76c^2k^4a+160c^2k^6b)}{32(a+4ck^2)(a+4bk^2)}. 
$$
In particular, when $a=0$ then 
$$ 
S(u,v) = -5ck^4/16 
$$ 
and when $b=0$ we obtain 
$$
S(u,v) = -\frac{k^2(7a^2+56ack^2+76c^2k^4)}{32(a+4ck^2)}. 
$$
In either case the sectional curvature is negative for any $k \neq 0$. 
All these examples work on $\mathbb{T}^n$ as well, if $x$ and $y$ denote the first 
two variables of the coordinate system.
\end{proof}

%%%%%%%%%%%%%%%%%%%%%%%%%%%%%%%%%%
\subsection{The $L^2$-metric on $\Diff(\mathbb{T}^n)$: the  Burgers equation} 
\label{subsec:burgers} 

If $b = c = 0$ then the formula \eqref{abcmetric} reduces to the $L^2$ inner product 
and the corresponding geodesic equation \eqref{utBuu} 
is the multi-dimensional Burgers equation 
\begin{equation*}
u_t + \nabla_uu + \diver(u) u 
+ 
\tfrac{1}{2} \nabla \langle u, u \rangle = 0, 
\end{equation*} 
also called the template matching equation; see \cite{HMA}. 
However, in contrast with the one-dimensional case, the curvature of 
the right-invariant $L^2$-metric assumes both signs when $n\ge 2$. 
For simplicity we only prove the result for $n=2$. 

\begin{proposition} \label{multiburgerscurvature}
The sectional curvature of $\Diff(\mathbb{T}^2)$ equipped with 
the right-invariant metric \eqref{abcmetric} with $b=c=0$  
is given by the formula \eqref{1dburgerscurvature} for any 
$u = f(x) \, \frac{\partial}{\partial x}$ 
and 
$v=g(x) \, \frac{\partial}{\partial x}$, in which case $S(u,v)\ge 0$. 
On the other hand, if 
$u=\sin{(2\pi x)} \, \frac{\partial}{\partial x}$ 
and 
$v = \sin^2{(2\pi x)} \, \frac{\partial}{\partial y}$ 
then $S(u,v) <0$. 
\end{proposition}
\begin{proof}
In this case the operator $\ad^\ast_v u$ defined by \eqref{opB} has the form 
$$
\ad^\ast_v u 
= 
u \diver{v} + (\iota_vdu^\flat)^{\sharp} 
+  
\grad \langle u, v \rangle.
$$
Thus, when $u=f(x) \frac{\partial}{\partial x}$ 
and 
$v=g(x) \tfrac{\partial}{\partial x}$ 
we get 
$\ad^\ast_v u 
= 
\big( 2g_x f + g f_x \big) \, \tfrac{\partial}{\partial x}$ 
which is the same formula as in the one-dimensional case and the first part of 
the proposition follows. 

Furthermore, if $w = g(x) \tfrac{\partial}{\partial y}$ then we compute 
\begin{align*}
\ad^\ast_w u 
&= 0, \\
\ad^\ast_u u 
&= 3f f_x  \tfrac{\partial}{\partial x}, \\
\ad^\ast_uw 
&= \big( f_x g 
+ 
f g_x \big)  \tfrac{\partial}{\partial y}, \\
\ad^\ast_ww
&= g g_x  \tfrac{\partial}{\partial x}.
\end{align*}
Combining these formulas with 
$\ad_uw = -f g_x \frac{\partial}{\partial y}$ 
in \eqref{generalcurvature}, we get 
$$ 
S(u,w) 
= 
a \int_0^1 \Big( \tfrac{1}{4} f_x^2 g^2 - 2 f f_x g g_x \Big) \, dx. 
$$ 
Taking $f(x)=\sin{(2\pi x)}$ and $g(x)=\sin^2{(2\pi x)}$, we find $S(u,w) = -15\pi^2/16$.
\end{proof}

A similar consideration in the general case can be summarized as the following statement.

\begin{theorem}\label{thm:multiburgerscurvature}
The sectional curvature of $\Diff(\mathbb{T}^n)$ equipped with 
the right-invariant $L^2$ metric (i.e. the $a$-$b$-$c$ metric \eqref{abcmetric} with $a=1$ and $b=c=0$)  
assumes both signs.
\end{theorem}

%%%%%%%%%%%%%%%%%%%%%%%%%%%%%%%%%%%%%%%%%%%%%%
\subsection{The homogeneous $\dot H^1$-metric on $\Diff_{\mu,\text{ex}}(\mathbb{T}^2)$} 
\label{subsec:abzero}

Consider next the Lie group  $\Diff_{\mu,\text{ex}}(\mathbb{T}^2)$ of exact volumorphisms 
of the flat torus $\mathbb{T}^2$ which consists of symplectic diffeomorphisms preserving 
the center of mass. 
Its Lie algebra consists of Hamiltonian vector fields $u= \sgrad{f}$ with $f\in C^\infty(\mathbb{T}^2)$. 
Following Arnold \cite{A} we calculate the sectional curvature of the metric
\begin{align}\label{cmetric}
\llangle u, v \rrangle 
= c\int_{\mathbb{T}^2} 
\langle du^{\flat}, du^{\flat}\rangle \, d\mu,
%du^\flat \wedge * dv^\flat, 
\end{align}
on this group. It turns out to be just as convenient to work with 
a more general right-invariant metric given at the identity by 
\begin{equation} \label{Lambda-metric} 
\llangle u,v\rrangle 
= 
\llangle \sgrad{f}, \sgrad{g} \rrangle 
= 
\int_{\mathbb{T}^2} f \Lambda g \, d\mu, 
\end{equation} 
where the (positive-definite, symmetric) operator defining the inner product 
is given by the formula 
$\Lambda = \lambda(\Laplacian)$ 
for some function 
$\lambda\colon \mathbb{R}^+\to\mathbb{R}^+$.\footnote{More precisely, 
if a function $f$ is written in an eigenbasis of the positive-definite 
Laplacian as 
$f=\sum_k a_k \phi_k$ with $\Laplacian \phi_k = \gamma_k \phi_k$, 
then 
$\Lambda f = \sum_k a_k \lambda(\gamma_k) \phi_k$.} 
For a vector $p \in \mathbb{R}^2$, we will write $F(p) = \lambda(\lvert p\rvert^2)$ 
for convenience. The metric (\ref{cmetric}) corresponds to $ \lambda(z)=cz^2$ 
and $\Lambda = c\Laplacian^2$.

\begin{theorem} \label{prop:lambda}     
Suppose $f(x,y) = \cos{(jx+ky)}$ and $g(x,y) = \cos{(lx+my)}$, where $j,k,l,m$ 
are integer multiples of $2\pi$. Set $p=(j,k)$ and $q=(l,m)$ and let
$u=\sgrad{f}$ and $v=\sgrad{g}$.
%$\Laplacian f = \lvert p \rvert^2 f$ and $\Laplacian g = \lvert q \rvert^2 g$ 
Then\footnote{Note that for $F(p)=\lvert p\rvert^2$ formula \eqref{curvaturevolumorphism} 
reproduces \eqref{arnoldeq}, up to a rescaling factor.}
\begin{multline}\label{curvaturevolumorphism}
S(u,v) 
= 
\frac{\lvert p\wedge q\rvert^2}{8} \bigg\{ 
\frac{1}{4}\Big( F(p)-F(q)\Big)^2 \Big( \frac{1}{F(p+q)} 
+ 
\frac{1}{F(p-q)}\Big)    \\ 
- 
\frac{3}{4} \Big( F(p+q)+F(p-q)\Big) + F(p) + F(q) \bigg\} 
\end{multline}
where $p\wedge q= jm - kl$.
\end{theorem}
\begin{proof}
%\begin{equation}\label{generalcurvature}
%\begin{split}
%\llangle R(u,v)v,u\rrangle 
%&= \frac{1}{4} \lVert B(u,v)+B(v,u)\rVert^2 
%- 
%\llangle B(u,u), B(v,v)\rrangle 
%\\ 
%&\qquad\qquad- \frac{3}{4} \lVert [u,v]\rVert^2 
%+ 
%\frac{1}{2} \llangle [u,v], B(u,v)-B(v,u)\rrangle.
%\end{split}
%\end{equation}
%
Recall that $\sgrad{f} = - f_y \tfrac{\partial}{\partial x} + f_x \tfrac{\partial}{\partial y}$ 
so that if $\{f,g\} = f_xg_y-f_yg_x$ denotes the Poisson bracket, then 
$[ \sgrad{f}, \sgrad{g} ] = \sgrad\{f,g\}$. 
Given any smooth functions $f$, $g$ and $h$ on the torus let 
$u=\sgrad{f}$, $v=\sgrad{g}$ and $w=\sgrad{h}$ be the corresponding skew-gradients. 
Integration by parts gives 
\begin{align*} 
\llangle \ad^\ast_v u, w\rrangle 
&= 
-\int_{\mathbb{T}^2} ( \Lambda f) \{ g,h\} \, d\mu    \\
&= 
\int_{\mathbb{T}^2} h \{ g, \Lambda f\} \, d\mu    
= 
\llangle w, \sgrad \Lambda^{-1} \{ g,\Lambda f \} \rrangle, 
\end{align*} 
from which we deduce that $\ad^\ast_v u = \sgrad \Lambda^{-1} \{g, \Lambda f\}$.

We can assume that 
$\lvert p\rvert$, $\lvert q\rvert$, $\lvert p-q\rvert$, and $\lvert p+q\rvert$ 
are all nonzero.
Furthermore, since 
$\Laplacian f = \lvert p\rvert^2 f$ and $\Laplacian g = \lvert q\rvert^2 g$ 
we find that 
\begin{align*} 
\ad^\ast_uu &= \ad^\ast_vv = 0, \\ 
\ad^\ast_v u &= -F(p) \sgrad \Lambda^{-1} \theta   
\quad 
\mathrm{and} 
\quad 
\ad^\ast_u v = F(q) \sgrad \Lambda^{-1}\theta, 
\end{align*} 
where 
\begin{align*}
\theta(x,y) 
= 
\{f,g\}(x,y)  
&= 
\frac{1}{2} (jm-kl)\Big( \cos{\big( (j-l)x+(k-m)y \big)}  \\ 
&\qquad\qquad\qquad\quad - 
\cos{\big( (j+l)x+(k+m)y\big) }\Big). 
\end{align*} 

Clearly, we have 
%Introducing functions $\varphi$ and $\psi$ such that 
$\theta=\frac{1}{2} (p\wedge q) (\varphi - \psi)$ 
%we compute 
where $\varphi$ and $\psi$ are eigenfunctions satisfying
$\Lambda\varphi = F(p-q)\varphi$ 
and 
$\Lambda \psi = F(p+q) \psi$. 
Combining the above formulas we obtain 
\begin{align*}
\ad_uv &= -\sgrad{\theta} = 
-\frac{p\wedge q}{2} \sgrad (\varphi - \psi), \\ 
\ad^\ast_v u &= -
\frac{(p\wedge q)F(p)}{2} \sgrad 
\left( \frac{\varphi}{F(p-q)} - \frac{\psi}{F(p+q)}\right),\\
\ad^\ast_u v &= 
\frac{(p\wedge q)F(q)}{2} \sgrad 
\left( \frac{\varphi}{F(p-q)} - \frac{\psi}{F(p+q)}\right).
\end{align*} 
Since, on the other hand, we have 
\begin{align*}
&\llangle \sgrad \varphi, \sgrad \varphi\rrangle 
= 
\tfrac{1}{2} F(p-q),  
\quad 
\llangle \sgrad \varphi, \sgrad \psi\rrangle 
= 0  \\
&\mathrm{and} \quad 
\llangle \sgrad \psi, \sgrad \psi\rrangle 
= \tfrac{1}{2} F(p+q),
\end{align*}
substituting into the sectional curvature formula \eqref{generalcurvature} 
yields \eqref{curvaturevolumorphism}.
\end{proof}

\begin{corollary}
In the particular case when the metric is given by \eqref{cmetric} 
the sectional curvature 
can assume both signs depending on $p$ and $q$.
\end{corollary}

Indeed, for the metric \eqref{cmetric} 
we have 
$F(p)=\lambda(\lvert p\rvert^2) = c\lvert p\rvert^4$ 
and a straightforward computation gives 
$S(u,v)>0$ when $p=(10\pi,0)$ and $q=(8\pi, 2\pi)$ 
while 
$S(u,v) <0$ when $p=(2\pi,0)$ and $q=(0,2\pi)$.

%%%%%%%%%%%%%%%%%%%%%%
%%%%%%%%%%%%%%%%%%%%%%%%%%

%%%%%%%%%%%%%%%%%%%%%%%%%%%

\appendix 
\section{Semi-direct products and curvature formulas} 
\label{semi} 
\nequation 
\numberwithin{equation}{section} 

Formula (\ref{abcurvature}) for the sectional curvature of the $H^1$ metric 
derived in Lemma \ref{1Dcurvprop} of Section \ref{circlediffeo} 
resembles the formula for the curvature
of a Riemannian submanifold $\mathcal{N}$ isometrically immersed in an 
ambient manifold $\mathcal{M}$, given by \eqref{gausscodazzi}.
Observe however that 
the sign of \eqref{gausscodazzi} is the opposite
of the sign for \eqref{abcurvature}. The next proposition shows how to 
rearrange \eqref{abcurvature} to make the analogy work. 

\begin{proposition} \label{curvaturerewrite}
The sectional curvature \eqref{abcurvature} of $\Diff(S^1)$ with the $H^1$ metric \eqref{abmetric} 
can be rewritten in either of the two following forms 
\begin{align} \label{mysteryembedding}
S(u,v) 
=& 
-\tfrac{a}{2} \int_{S^1} ( uv_x - vu_x )^2 \, dx 
+ 
\tfrac{a}{b} \llangle T(u,u), T(v,v)\rrangle - \tfrac{a}{b} \llangle T(u,v), T(u,v)\rrangle 
	\\  \label{semidirectembedding}
S(u,v) 
= &\; 
\tfrac{1}{a} \int_{S^1} \big(a(uv_x -vu_x ) 
+ 
\tfrac{b}{2}( v_x u_{xx} - u_x v_{xx} ) \big)^2 \, dx 
	\\ \nonumber
&+ 
\tfrac{b}{a} \llangle Q(u,u), Q(v,v)\rrangle - \tfrac{b}{a} \llangle Q(u,v), Q(u,v)\rrangle,
\end{align}
where the bilinear maps $T$ and $Q$ are defined by
\begin{align} \label{TQdef} 
T(u,v) = A^{-1}\left(auv+\tfrac{b}{2}u_xv_x\right), 
\qquad 
Q(u,v) = -\partial_x^2A^{-1}\left(auv+ \tfrac{b}{2}u_xv_x\right).
\end{align}
\end{proposition}
\begin{proof}
Note that $\Gamma(u,v)$ defined by \eqref{1DGammadef} is related to $Q$ and $T$ by 
$\Gamma(u,v) = \partial_x T(u,v)$ and $Q(u,v) = -\partial_x \Gamma(u,v)$. 

We prove \eqref{mysteryembedding} first. 
Let $u$ and $v$ be any vector fields and set $q=auv+\frac{b}{2} u_xv_x$. 
Then $T(u,v) = A^{-1}q$, so that
\begin{align*} 
\llangle \Gamma(u,v),\Gamma(u,v)\rrangle 
&= 
\int_{S^1} \partial_x T(u,v) A \partial_x T(u,v) \, dx 
= 
-\int_{S^1} q A^{-1} \partial_x^2 q \, dx.
\end{align*}
Now using the identity $A^{-1}\partial_x^2 = -\frac{1}{b} + \frac{a}{b} A^{-1}$, 
we find 
$$ 
\llangle \Gamma(u,v), \Gamma(u,v)\rrangle 
= 
\frac{1}{b} \int_{S^1} \left(auv+ \tfrac{b}{2}u_xv_x\right)^2 \, dx 
- 
\frac{a}{b} \llangle T(u,v), T(u,v)\rrangle. 
$$ 
Using the same trick also on the other term in \eqref{abcurvature} gives 
\eqref{mysteryembedding}.
A similar technique reduces \eqref{semidirectembedding} to \eqref{abcurvature}.
\end{proof} 

The formulas \eqref{mysteryembedding} and \eqref{semidirectembedding} suggest that 
$\Diff(S^1)$ with the $a$-$b$ metric \eqref{abmetric} can be realized as 
an isometrically immersed submanifold of some simpler manifold. 
For \eqref{semidirectembedding} such an immersion is described in the following theorem. 

\begin{theorem} \label{semidirecttheorem}
Let $\Diff(S^1)\ltimes C^{\infty}(S^1)$ denote the semidirect product of 
the diffeomorphism group $\Diff(S^1)$ with $C^{\infty}(S^1)$ 
with group structure given by 
$$ 
(\eta, F)\cdot (\xi, G) = (\eta\circ\xi, F\circ\xi + G). 
$$
Define a right-invariant Riemannian metric on this group by the $L^2$ inner product 
at the identity 
\begin{equation} \label{semidirectL2}
\llangle (u, f), (v, g)\rrangle = \int_{S^1} \big( auv + bf g \big) dx.
\end{equation} 
Let $\Upsilon \colon \Diff(S^1)\to \Diff(S^1)\ltimes C^{\infty}(S^1)$ denote the map 
$\Upsilon(\eta) = (\eta, \ln{\eta_x })$. 
Then $\Upsilon$ is an embedding and a group homomorphism 
and the right-invariant metric induced on $\Diff(S^1)$ by \eqref{semidirectL2} 
is the $a$-$b$ metric \eqref{abmetric}. 
Furthermore, the curvature formula \eqref{semidirectembedding} 
is the Gauss-Codazzi formula for the embedding $\Upsilon$. 
\end{theorem} 

\begin{proof}
The Lie algebra for $G=\Diff(S^1)\ltimes C^{\infty}(S^1)$ is the semidirect product 
$\mathfrak{g} = \Vect(S^1)\ltimes C^{\infty}(S^1)$, 
which has an interpretation as the space of first-order differential operators 
$vD + f$ for $v\in \Vect$ and $f\in C^{\infty}$. 
For details on the geometry of the semidirect product see e.g., \cite{MRW, vizman}. 
For this semi-direct product  the adjoint operator has the form 
$$ 
\ad_{(u,f)}(v,g) = (-uv_x +u_x v, vf_x -ug_x ), 
$$
which implies that the Arnold operator \eqref{opB} is 
\begin{equation}\label{Bopsemidirect}
\ad_{(u,f)}^{\ast}(v,g) 
= 
\left( 2u_x v+uv_x +\tfrac{b}{a} gf_x , gu_x +g_x u \right).
\end{equation} 
The general curvature formula \eqref{generalcurvature} implies, after some simplifications, 
that the curvature of the semidirect product $\Diff(S^1)\ltimes C^{\infty}(S^1)$ 
can be written as
\begin{multline} \label{semidirectcurvature}
\llangle \bar{R}((u,f), (v,g))(v,g), (u,f)\rrangle 
= 
a \int_{S^1} \big( uv_x -vu_x + \tfrac{b}{2a} (gf_x -fg_x )\big)^2 \, dx \\
 + \tfrac{b}{4} \int_{S^1} (gu_x -fv_x ) (gu_x -fv_x + 8vf_x -8ug_x ) \, dx.
\end{multline}
Observe that  $\Upsilon$ is a group homomorphism: 
$$ 
\Upsilon(\eta\circ\xi) 
= 
\big(\eta\circ\xi, \ln(\eta_x \circ\xi) + \ln( \xi_x )\big) 
= 
\Upsilon(\eta)\cdot \Upsilon(\xi) 
$$
using the chain rule. 
Smoothness of this map can be proved using Sobolev $H^s$ topology. 
For our purposes it is enough to note that $\Upsilon$ is formally an immersion 
since at the identity we have 
$$ 
D\Upsilon_e(u) = ( u, u_x), 
$$ 
which is obviously injective.\footnote{In the interpretation via differential operators, we consider operators of the form $uD+u_x=Du$. 
The homomorphism $ D\Upsilon: uD\mapsto Du$ is evidently a homomorphism of 
Lie algebras $\Vect\to \Vect\ltimes C^{\infty}$.}
Hence the induced metric on $\Diff(S^1)$ is given at the identity by 
$$ 
\llangle u, v\rrangle 
= 
\llangle D\Upsilon(u), D\Upsilon(v)\rrangle 
= 
\int_{S^1} \big( a u v + bu_x v_x \big) dx, 
$$
which is precisely the $a$-$b$ metric \eqref{abmetric}.

The orthogonal complement of the image of $D\Upsilon_e$ in the metric \eqref{semidirectL2} 
consists of those vectors of the form 
$\big(\tfrac{b}{a} h_x, h \big)$ 
for some function $h\colon S^1\to \mathbb{R}$. 
We can now compute the second fundamental form \eqref{secondfundy} of the embedding: 
because of the right-invariance, we have 
$\llangle \Pi(u,u), w\rrangle = \llangle \ad^\ast_uu, w\rrangle$ 
for 
$w \in (\text{Im}\; D\Upsilon_e)^\perp$ 
where $\ad^\ast$ is the operator given in \eqref{Bopsemidirect}. 
We then have 
$$ 
\Pi(u,u) 
= 
\left(\tfrac{b}{a}  h_x , h\right) \quad \text{where}\quad h 
= 
-A^{-1} \partial_x^2 \left(au^2 + \tfrac{b}{2}u_x^2\right). 
$$
Polarization yields 
$\Pi(u,v) = \left( \tfrac{b}{a} Q(u,v)_x , Q(u,v)\right)$, 
where $Q$ is given by (\ref{TQdef}). 
Substituting \eqref{semidirectcurvature} when $f=u'$ and $g=v'$, 
together with the second fundamental form, into \eqref{gausscodazzi} reproduces 
\eqref{semidirectembedding}.
\end{proof}

We conjecture that there is another embedding-homomorphism of $\Diff(S^1)$ into an infinite-dimensional Lie group with right-invariant metric which reproduces formula \eqref{mysteryembedding}, but we do not know what it might be.
%
%\begin{remark}\upshape
%  The geodesic equation $u_t = -\ad_u^{\ast}u$ on the semidirect product $\Diff(S^1)\ltimes C^{\infty}(S^1)$ equipped with the metric (\ref{semidirectL2}) is
%  $$u_t = -3uu_x - \frac{b}{a}ff_x, \qquad
%  f_t = -(uf)_x.$$ 
%  This bears a superficial resemblance to the barotropic equations of compressible fluid mechanics, although one cannot eliminate the factor $3$ by a rescaling. 
%\end{remark}

%%%%%%%%%%%%%%%%%%%%%%%%
%%%%%%%%%%%%%%%%%%%%%%%%
%%%%%%%%%%%%%%%%%%%%%%%%

\bigskip
{\bf Acknowledgements.}
B.K. was partially supported by an NSERC
research grant. 
J.L. acknowledges support from the EPSRC, UK. 
G.M. was supported in part by The James D. Wolfensohn Fund and Friends of the Institute for 
Advanced Study grants.
S.C.P. was partially supported by an NSF grant.

\bibliographystyle{plain}

\bibliography{is}

\end{document}